\documentclass[aos,seceqn,MSNbibl,citesort,dvips]{arximspdf}

\doi{10.1214/10-AOS861}
\volume{39}
\issue{3}
\pubyear{2011}
\firstpage{1526}
\lastpage{1550}

\makeatletter

\newtheorem{theorem}{Theorem}[section]
\newtheorem{lemma}{Lemma}[section]
\newproclaim{remark}{Remark}

\newcommand{\boldmathxione}{\bolds\xi_{1}}
\newcommand{\boldmathxitwo}{\bolds\xi_{2}}
\newcommand{\boldeta}{\bolds\eta}
\newcommand{\boldTheta}{\bolds\Theta}
\newcommand{\boldxi}{\bolds\xi}

\makeatother

\begin{document}
\begin{frontmatter}

\title{Uniform moment bounds of Fisher's information with applications to time series}
\runtitle{Uniform moment bounds}

\begin{aug}
\author[A]{\fnms{Ngai Hang} \snm{Chan}\thanksref{t1}\ead[label=e1]{nhchan@sta.cuhk.edu.hk}} and
\author[B]{\fnms{Ching-Kang} \snm{Ing}\corref{}\thanksref{t2}\ead[label=e2]{cking@stat.sinica.edu.tw}}
\runauthor{N. H. Chan and C.-K. Ing}
\affiliation{Chinese University of Hong Kong and Academia Sinica}
\address[A]{Department of Statistics\\
Room 118, Lady Shaw Building\\
Chinese University of Hong Kong\\
Shatin, New Territories\\
Hong Kong \\
\printead{e1}}
\address[B]{Institute of Statistical Science\\
Academia Sinica\\
Taipei 115\\
Taiwan, ROC\\
\printead{e2}}
\end{aug}

\thankstext{t1}{Supported in part by the General Research Fund Nos.
400408 and 400410 and
the Collaborative Research Fund No. CityU8/CRF/09, all from the
Research Grants Council of Hong Kong.}

\thankstext{t2}{Supported in part by the National Science Council of
Taiwan under Grant NSC 94-2118-M-001-013.}

\received{\smonth{5} \syear{2010}}
\revised{\smonth{8} \syear{2010}}

%
\begin{abstract}
In this paper, a uniform (over some parameter space) moment bound for
the inverse of Fisher's information matrix is established. This result
is then applied to develop moment bounds for the normalized least
squares estimate in (nonlinear) stochastic regression models. The
usefulness of these results is illustrated using time series models. In
particular, an asymptotic expression for the mean squared prediction
error of the least squares predictor in autoregressive moving average
models is obtained. This asymptotic expression provides a solid
theoretical foundation for some model selection criteria.
\end{abstract}

%
\begin{keyword}[class=AMS]
\kwd[Primary ]{62J02}
\kwd[; secondary ]{62M10}
\kwd{62F12}
\kwd{60F25}.
\end{keyword}
\begin{keyword}
\kwd{Fisher's information matrix}
\kwd{least squares estimates}
\kwd{mean squared prediction errors}
\kwd{stochastic regression models}
\kwd{uniform moment bounds}.
\end{keyword}

\end{frontmatter}

\section{Introduction}\label{intro}

Moment inequalities and moment bounds have long been vibrant topics in
modern probability and statistics. The celebrated inequalities of Burkholder
\cite{3} and Doob~\cite{5} offer exemplary illustrations of the
importance of
moment inequalities. Using moment bounds, the
order of magnitude of the spectral norm of the inverse of the Fisher's
information
matrix can be quantified and consistency and efficiency of least
squares estimates
of stochastic regression and adaptive control can be established; see,
for example, the seminal work of
Lai and Wei~\cite{16} and the succinct review of Lai and Ying \cite
{15}. In this paper, a uniform
(over some parameter space) moment bound for
the inverse of the Fisher's information matrix is established. This bound
is used to investigate the moment properties of least squares
estimates and the mean squared prediction error (MSPE) for time series models.

To appreciate the significance of uniform moment bounds, consider
the stochastic regression model
%
%
\begin{equation}\label{equ1.1}
y_{t}=g_{t}(\theta_{0})+\varepsilon_{t}, \qquad t=1, \ldots, n,
\end{equation}
where $g_{t}(\cdot)$ is a random function,
$\theta_{0}$ is an unknown parameter and $\{\varepsilon_{t}\}$
is a martingale difference sequence. There are two important problems
related to this model.

The first one concerns the mean squared error prediction. In practice,
the unknown parameter $\theta_{0}$ is usually estimated by the least
squares estimate $\hat{\theta}_{n}$,
which minimizes $S_{n}(\theta)=\sum_{t=1}^{n}(y_{t}-g_{t}(\theta))^{2}$.
Although the (strong) law of large numbers (LLN) and the central limit
theorem (CLT)
of $\hat{\theta}_{n}$ were established under certain assumptions on
$g_{t}(\cdot)$ and $\varepsilon_{t}$
(see among others, Lai~\cite{14} and Skouras~\cite{19}),\vspace*{2pt} relatively
little is known about
the moment convergence of $\hat{\theta}_{n}$. Moment convergence of
$\hat{\theta}_{n}$ offers important insight
in the pursuit of the mean squared prediction problem.
To see this, suppose that $n^{1/2}(\hat{\theta}_{n}-\theta_{0})$ is
asymptotically normal with mean zero
and variance $\eta> 0$. Then an immediate question is to pursue
%
%
\begin{equation}\label{equ1.2}
\mathrm{E}|n^{1/2}(\hat{\theta}_{n}-\theta_{0})|^{q}=O(1),\qquad q \geq1.
\end{equation}
In particular, if (\ref{equ1.2}) holds for some $q>2$, then $\{n(\hat
{\theta}_{n}-\theta_{0})^{2}\}$ is uniformly
integrable and consequently,
$\lim_{n \to\infty} n {\mathrm E}(\hat{\theta}_{n}-\theta
_{0})^{2}=\eta$.
This result can be applied to develop an asymptotic expression for the
mean squared error of $\hat{\theta}_{n}$ as
\[
{\mathrm E}(\hat{\theta}_{n}-\theta_{0})^{2}= \frac{\eta}{n}+o(n^{-1})
\]
from which asymptotic properties of the MSPE
of the least squares predictor $g_{n+1}(\hat{\theta}_{n})$ of $y_{n+1}$,
$E(y_{n+1}-g_{n+1}(\hat{\theta}_{n}))^{2}$,
can be established; see Sections~\ref{sec2} and~\ref{sec3} for further details.

To establish (\ref{equ1.2}), consider the Fisher's information
number,\break
$n^{-1}\sum_{t=1}^{n}(g^{\prime}_{t}(\theta))^{2}$ of (\ref{equ1.1}),
where $g^{\prime}_{t}(\theta)=d g_{t}(\theta)/d \theta$.
As will be shown in Section~\ref{sec2}, it turns out that the uniform negative
moment bound for\break
$n^{-1}\sum_{t=1}^{n}(g^{\prime}_{t}(\theta))^{2}$, that is, for any $q
\geq1$,
%
%
\begin{equation}\label{equ1.3}
{\mathrm E}\Biggl\{ \sup_{\theta\in
B_{\delta_{1}}(\theta_{0})}\Biggl(n^{-1}\sum_{t=1}^{n}(g^{\prime}_{t}(\theta
))^{2}\Biggr)^{-q}\Biggr\}=O(1)
\end{equation}
plays a crucial role in proving (\ref{equ1.2}),
where $B_{\delta_{1}}(\theta_{0})=\{\theta\dvtx |\theta-\theta
_{0}|<\delta
_{1}\}$ for some $\delta_{1}>0$.

A second but equally important problem in stochastic regression
concerns model selection. To understand how
the uniform moment bound is related to this issue,
consider the case when $g_{t}(\cdot)$ in (\ref{equ1.1}) contains $k>1$
unknown parameters $\bolds\theta_{0} \in R^{k}$.
A multiparameter generalization of (\ref{equ1.3}) becomes: for any $q
\geq1$,
%
%
\begin{equation}\label{equ1.4}
{\mathrm E}\Biggl\{ \sup_{\bolds\theta\in
B_{\delta_{1}}(\bolds\theta_{0})}
\lambda^{-q}_{\min}\Biggl(n^{-1}\sum_{t=1}^{n} \nabla g_{t}(\bolds\theta)
(\nabla g_{t} (\bolds\theta))^{\mathrm T}\Biggr)\Biggr\}= O(1),
\end{equation}
where $\lambda_{\min}(L)$ denotes the minimum eigenvalue of the
matrix $L$
and $\nabla g_{t} (\bolds\theta)$ denotes the gradient vector of
$g_{t}(\bolds\theta)$.
In particular, when $g_{t}(\bolds\theta)=g_{t}(\theta_{1},
\ldots,\break\theta
_{k})=\theta_{1}y_{t-1}+\cdots+\theta_{k}y_{t-k}$ in
(\ref{equ1.1}), that is, when $y_t$ is an autoregressive (AR) model of
order $k$, (\ref{equ1.4}) reduces to
%
%
\begin{equation}\label{equ1.5}
{\mathrm E}\Biggl\{
\lambda^{-q}_{\min}\Biggl(n^{-1}\sum_{t=1}^{n}
\mathbf{y}_{t-1}(k) \mathbf{y}^{\mathrm T}_{t-1} (k)\Biggr)\Biggr\}= O(1),
\end{equation}
where $\mathbf{y}_{t}(k)=(y_{t}, \ldots, y_{t-k+1})^{\mathrm T}$.
By imposing a Lipschitz type condition on the distribution function of
$\varepsilon_{t}$
and a stationarity condition on $g_{t}(\cdot)$,
Findley and Wei~\cite{7} established (\ref{equ1.5}), thereby providing
a rigorous mathematical derivation of the AIC model selection criterion
for weakly stationary AR processes. However, the proof of (\ref
{equ1.4}) for a general stochastic
regression model is much more involved than (\ref{equ1.5})
due to the presence of an ``extra'' supremum, which is taken over an
uncountable set inside the expectation. As a
consequence, similar to the AR case, knowledge about negative uniform
moment bounds of the Fisher's
information matrix (\ref{equ1.4}) constitutes an indispensable tool for
the model selection problem.

The rest of this paper is organized as follows.
In Section~\ref{sec2}, we first show in Theorem~\ref{theo2.1} that (\ref{equ1.4}) holds
under more general situations
where $B_{\delta_{1}}(\bolds\theta_{0})$ is replaced by a bounded subset
$\boldTheta$ of $R^{k}$
and $\nabla g_{t}(\bolds\theta)$ is replaced by a vector-valued random
function $\mathbf{f}_{t}(\bolds\theta),
\bolds\theta\in\boldTheta$, satisfying certain assumptions.
We then apply Theorem~\ref{theo2.1} to establish the moment convergence of least
squares estimates
in (nonlinear) stochastic regression models; see Theorem~\ref{theo2.2}.
Section~\ref{sec3} focuses on the applications of Theorems~\ref{theo2.1} and~\ref{theo2.2} to
autoregressive moving average (ARMA) models.
In particular, the moment convergence of the least squares estimates
and an asymptotic expression (up to terms of order $n^{-1}$)
for the MSPE of the least squares predictor for ARMA models are established.
To facilitate the presentation,
technical results of Sections~\ref{sec2} and~\ref{sec3} are deferred to Appendices~\ref{appA} and
\ref{appB}, respectively.

\section{Uniform bounds on negative moments}\label{sec2}
Let $(\Omega, \mathcal{F}, P)$ be a probability space and
$\{\mathcal{F}_{t}\}$
be an increasing sequence of $\sigma$-fields on $(\Omega, \mathcal{F},
P)$. Let
$\mathbf{f}_{t}(\bolds\theta), t=1, \ldots, n$,
be $r$-dimensional $\mathcal{F}_{t}$-measurable
random functions of a parameter vector $\bolds\theta=
(\theta_{1}, \ldots, \theta_{k})^{\mathrm T} \in\boldTheta\subset R^{k}$.
In the first half this section, we provide sufficient conditions
under which the minimum eigenvalue of
the normalized matrix
$n^{-1} \sum_{t=1}^{n}\mathbf{f}_{t}(\bolds\theta)
\mathbf{f}_{t}^{\mathrm T}(\bolds\theta)$,
$\lambda_{\min}(
n^{-1} \sum_{t=1}^{n}\mathbf{f}_{t}(\bolds\theta)
\mathbf{f}_{t}^{\mathrm T}(\bolds\theta)
)$, satisfies the following uniform moment bound:
%
%
\begin{equation}\label{equ2.1}
\mathrm{E}\Biggl\{ \sup_{\bolds\theta
\in\boldTheta} \lambda_{\min}^{-q}\Biggl(n^{-1}
\sum_{t=1}^{n}\mathbf{f}_{t} (\bolds\theta)
\mathbf{f}_{t}^{\mathrm T}(\bolds\theta)\Biggr)\Biggr\} = O(1)\qquad \mbox{for any }
q \geq1.
\end{equation}
This uniform negative moment bound is applied to
investigate the moment properties of least squares
estimates in the second half of this section. To begin, assume the
following conditions:

\begin{longlist}[(C4)]
\item[(C1)]
$\mathbf{f}_{t}(\bolds\theta)$ is continuous
on $\boldTheta$ and $\boldTheta$
is a bounded subset of $R^{k}$;

\item[(C2)]
there exist positive integer $d$ and positive numbers $\delta$,
$\alpha$ and $M$ such that for any $t > d$, any $0 < s_{2}
- s_{1} \leq\delta$, any $\bolds\theta\in\boldTheta$ and any
$\|\mathbf{a}\|=1$,
\[
P\bigl(s_{1} < \mathbf{a}^{\mathrm T} \mathbf{f}_{t}
(\bolds\theta) \leq s_{2}|
\mathcal{F}_{t-d}\bigr) \leq M (s_{2} - s_{1})^{\alpha} \qquad\mbox{a.s.},
\]
where $\|\mathbf{a}\|$ denotes the Euclidean norm of vector
$\mathbf{a} \in R^{r}$;

\item[(C3)] there exist $\tau>0$
and nonnegative random variables $B_{t}$ satisfying
$\sup_{t \geq1} \mathrm{E}(B_{t}) \leq C_{1}$ for some $C_{1}>0$
such that for all
$\boldmathxione, \boldmathxitwo\in
\boldTheta$ with $\|\boldmathxione
-\boldmathxitwo\| < \tau$,
\[
\|\mathbf{f}_{t}(\boldmathxione)-\mathbf{f}_{t}
(\boldmathxitwo)\| \leq B_{t}\|\boldmathxione
-\boldmathxitwo\| \qquad\mbox{a.s.};
\]

\item[(C4)]
there exists $C_{2}>0$ such that
$\sup_{t \geq1} \mathrm{E}(\sup_{\bolds\theta
\in\boldTheta}
\|\mathbf{f}_{t}(\bolds\theta)\|^{2}) \leq C_{2}$.
\end{longlist}

(C1) is a standard assumption for the regression function and its
gradient vector
in nonlinear regression; see, for example, Lai~\cite{14} and Robinson
and Hidalgo~\cite{17}.
(C2) says that given the information ($\sigma$-field)
whose time index is sufficiently smaller than the current time index $t$,
the conditional distribution of $\mathbf{a}^{\mathrm T} \mathbf
{f}_{t}(\bolds\theta)$
follows a local Lipschitz condition of order $\alpha$ for all
points $\bolds\theta\in\boldTheta$
and all directions $\mathbf{a}$ with $\|\mathbf{a}\|=1$.
In the special case when $\boldTheta$ contains only
one point, (C2) is related to Findley and Wei's~\cite{7}
\textit{uniform Lipschitz condition over all directions},
which is the key assumption used in deriving the AIC
for stationary AR models.
Since we need to deal with the supremum over a class
of inverses of minimum eigenvalues indexed by $\bolds\theta$,
a Lipschitz type condition over all points ($\bolds\theta$)
in all directions ($\mathbf{a}$) is required in this paper.
As will be seen in Section~\ref{sec3},
(C2) is flexible enough to encompass many time series applications.
Conditions like (C3) have been imposed
on the regression function by Andrews~\cite{2} and Skouras~\cite{19}
in proving
the uniform law of large numbers for random functions
associated with $S_{n}(\theta)$.
(C3) can be verified when $\mathbf{f}_{t}(\bolds\theta)$
is sufficiently smooth; see (\ref{equ3.26}) for more details.
(C4) imposes a mild moment condition on $\mathbf{f}_{t}(\bolds\theta)$
and appears to be satisfied in many practical situations.
Moreover, (C4) can be weakened to
$\sup_{t \geq1}\sup_{\bolds\theta\in\boldTheta}\mathrm{E}(\|
\mathbf
{f}_{t}(\bolds\theta)\|^{2})\leq C_{2}$
for some $C_{2}>0$ at the price of strengthening the conditions on
$B_{t}$ in (C3) to
$\sup_{t \geq1}\mathrm{E}(B^{2}_{t})\leq C_{1}$
for some $C_{1}>0$.

\begin{theorem}\label{theo2.1} Assume that \textup{(C1)}--\textup{(C4)} hold. Then
inequality (\ref{equ2.1}) is true.\vadjust{\goodbreak}
\end{theorem}
\begin{pf}
First, note that the measurability of
$\sup_{\bolds\theta
\in\boldTheta} \lambda_{\min}^{-q}(n^{-1}
\sum_{t=1}^{n}\mathbf{f}_{t} (\bolds\theta)\times
\mathbf{f}_{t}^{\mathrm T}(\bolds\theta))$
is ensured by the continuity
of $\mathbf{f}_{t}(\bolds\theta)$.
Define $n_{d}= \lfloor(n-d)/d \rfloor$,
where $\lfloor a \rfloor$
is the largest integer \mbox{$\leq$}$a$.
Then for $n$ large,
%
%
\begin{eqnarray}\label{equ2.2}
&&
n^{q} \lambda_{\min}^{-q} \Biggl(\sum_{t=1}^{n}
\mathbf{f}_{t}(\bolds\theta)\mathbf{f}^{\mathrm T}_{t}
(\bolds\theta)\Biggr)\nonumber\\
&&\qquad\leq n^{q} \Biggl\{\sum_{j=1}^{d}\lambda_{\min} \Biggl(\sum_{i=0}^{n_{d}-1}
\mathbf{f}_{(i+1)d+j}(\bolds\theta)\mathbf{f}^
{\mathrm T}_{(i+1)d+j}(\bolds\theta)\Biggr)\Biggr\}^{-q} \\
&&\qquad\leq
\{n/(n_{d}d)\}^{q} d^{-1} \sum_{j=1}^{d} n^{q}_{d}\lambda^{-q}_{\min}
\Biggl(\sum_{i=0}^{n_{d}-1}\mathbf{f}_{(i+1)d+j}(\bolds\theta)
\mathbf{f}^{\mathrm T}_{(i+1)d+j}(\bolds\theta)\Biggr),\nonumber
\end{eqnarray}
where the first inequality is ensured by the fact that for symmetric
matrices $E_{1}$ and~$E_{2}$,
$\lambda_{\min}(E_{1}+E_{2}) \geq\lambda_{\min}(E_{1})+\lambda
_{\min}(E_{2})$,
and the second one is ensured by the convexity of $x^{-q}, x >0$.
As a key step for achieving (\ref{equ2.1}),
we show, by making use of (C2)--(C4), in Appendix~\ref{appA} that
there exists a positive integer $m$,
depending only on $q, r, k$ and $\alpha$, such that
for all large $n$, all $0\leq l \leq n_{d}-m$ and all $1 \leq j \leq d$,
%
%
\begin{equation}\label{equ2.3}
\mathrm{E}\Biggl(\sup_{\bolds\theta\in\boldTheta}
\lambda^{-q}_{\min}
\Biggl(\sum_{i=l}^{l+m-1}\mathbf{f}_{(i+1)d+j}(\bolds\theta)
\mathbf{f}^{\mathrm T}_{(i+1)d+j}(\bolds\theta)\Biggr)\Biggr) \leq C_{3},
\end{equation}
where $C_{3}$ is some positive constant independent of $l$ and $j$.
Let $n_{d,m}=\lfloor n_{d}/m \rfloor$. Then, analogous to (\ref{equ2.2}),
\begin{eqnarray*}
&& n^{q}_{d}\lambda^{-q}_{\min}
\Biggl(\sum_{i=0}^{n_{d}-1}\mathbf{f}_{(i+1)d+j}(\bolds\theta)
\mathbf{f}^{\mathrm T}_{(i+1)d+j}(\bolds\theta)\Biggr) \\
&&\qquad\leq(n_{d}/n_{d, m})^{q}
n^{-1}_{d, m} \sum_{s=0}^{n_{d, m}-1}
\lambda^{-q}_{\min}
\Biggl(\sum_{i=0}^{m-1}\mathbf{f}_{(i+sm+1)d+j}(\bolds\theta)
\mathbf{f}^{\mathrm T}_{(i+sm+1)d+j}(\bolds\theta)\Biggr).
\end{eqnarray*}
Combining this fact with (\ref{equ2.2}) and (\ref{equ2.3}) yields
for $n$ large and for some positive number $C_{4}$,
\begin{eqnarray*}
&& n^{q} \mathrm{E}\Biggl\{\sup_{\bolds\theta
\in\boldTheta}
\lambda_{\min}^{-q} \Biggl(\sum_{t=1}^{n}
\mathbf{f}_{t}(\bolds\theta)\mathbf{f}^{\mathrm T}_{t}
(\bolds\theta)\Biggr)\Biggr\}\\
&&\qquad\leq\frac{n^{q}}{(n_{d, m}d)^{q}d} \\
&&\qquad\quad{}\times\sum_{j=1}^{d}
n^{-1}_{d, m}
\sum_{s=0}^{n_{d, m}-1}
\mathrm{E}\Biggl\{\sup_{\bolds\theta\in\boldTheta}
\lambda^{-q}_{\min}
\Biggl(\sum_{i=0}^{m-1}\mathbf{f}_{(i+sm+1)d+j}(\bolds\theta)
\mathbf{f}^{\mathrm T}_{(i+sm+1)d+j}(\bolds\theta)\Biggr)\Biggr\}\\
&&\qquad\leq C_{3}C_{4}m^{q}.
\end{eqnarray*}
Thus, (\ref{equ2.1}) follows.
\end{pf}

To see the extent of the usefulness of (\ref{equ2.1}),
consider a stochastic regression model of the form
%
%
\begin{equation}\label{equ2.4}
y_{t}=g_{t}(\bolds\theta_{0})+\varepsilon_{t},\qquad t=1, \ldots, n,
\end{equation}
where
$\{\varepsilon_{t}\}$
is a martingale difference sequence
with respective to $\{\mathcal{G}_{t}\}$, an increasing sequence of
$\sigma$-fields on $(\Omega, \mathcal{F}, P)$,
such that
%
%
\begin{equation}\label{equ2.5}
\sup_{t}\mathrm{E}(\varepsilon^{2}_{t}|\mathcal{G}_{t-1})< \infty\qquad
\mbox{a.s.},
\end{equation}
$g_{t}(\cdot)$
is a $\mathcal{G}_{t-1}$-measurable
random function on a compact set $\boldTheta_{1} \subset R^{k}$
and $\bolds\theta_{0} \in\boldTheta_{1}$
is unknown coefficient vector.
The least squares estimate
$\hat{\bolds\theta}_{n}$ of $\bolds\theta_{0}$
is obtained by
minimizing
%
%
\begin{equation}\label{equ2.6}
S_{n}(\bolds\theta)=\sum_{t=1}^{n}\bigl(y_{t}-g_{t}
(\bolds\theta)\bigr)^{2}
\end{equation}
over $\boldTheta_{1}$.
The next theorem provides a set of sufficient conditions under which
%
%
\begin{equation}\label{equ2.7}
\mathrm{E}\|n^{1/2}(\hat{\bolds\theta}_{n}-\bolds\theta_{0})\|
^{q} = O(1),\qquad
q \geq1.
\end{equation}
To state the result, denote the gradient vector and the Hessian
matrix of a smooth function $h\dvtx R^{k} \to R$
by $\nabla h (\xi_{1}, \ldots, \xi_{k})= (\partial h/\partial\xi_{1},
\ldots, \partial h/\partial\xi_{k})^{\mathrm T}$
and $\nabla^{2} h (\xi_{1}, \ldots, \xi_{k})=
(\partial^{2}h/\partial\xi_{i}\, \partial\xi_{j})_{1\leq i, j\leq k}$,
respectively.
For $\bolds\theta\in R^{k}$
and $\eta_{1}>0$,
define $B_{\eta_{1}}(\bolds\theta)=
\{\boldxi\dvtx\|\boldxi-
\bolds\theta\|< \eta_{1} \}$.
\begin{theorem}\label{theo2.2}
Consider the stochastic regression model (\ref{equ2.4})
in which $g_{t}(\cdot)$ is $\mathcal{G}_{t-1}$-measurable and
continuous on $\boldTheta_{1}$ and the
martingale difference sequence $\{\varepsilon_{t}\}$ satisfies
(\ref{equ2.5}).
Suppose that there exists $\delta_{1}>0$ such that
$B_{\delta_{1}}(\bolds\theta_{0}) \subset\boldTheta_{1}$ and
the gradient vector $\nabla g_{t}$ is continuously differentiable on
$B_{\delta_{1}}(\bolds\theta_{0})$. Moreover, assume
$\sup_{t}\mathrm{E}(|\varepsilon_{t}|^{\gamma}|\mathcal{G}_{t-1})<
C_{5}$ a.s. for some $\gamma>\max\{q, 2\}$ and $C_{5}>0$,
and the following conditions hold:
\begin{longlist}
\item
\textup{(C2)}--\textup{(C4)} hold for $\boldTheta=
B_{\delta_{1}}(\bolds\theta_{0})$,
$\mathbf{f}_{t}(\bolds\theta)=
\nabla g_{t} (\bolds\theta)$ and
$\mathcal{F}_{t}=\mathcal{G}_{t-1}$.
In addition, there exists $q_{1}>q$ such that
%
%
\begin{eqnarray}\label{equ2.8}
\max_{1\leq i, j \leq k}\mathrm{E}\Biggl(\sup_{\bolds\theta\in
B_{\delta_{1}}(\bolds\theta_{0})}
\Biggl|n^{-1/2}\sum_{t=1}^{n}\varepsilon_{t}(\nabla^{2} g_{t}
(\bolds\theta))_{i,j}\Biggr|^{q_{1}}\Biggr)&=& O(1),
\\
\label{equ2.9}
\max_{1\leq i, j \leq k, 1 \leq t \leq n}\mathrm{E}\Bigl(\sup_
{\bolds\theta\in B_{\delta_{1}}(\bolds\theta_{0})}
|(\nabla^{2} g_{t}(\bolds\theta))_{i,j}|^{4 q_{1}}\Bigr)&=&O(1),
\\
\label{equ2.10}
\max_{1 \leq t \leq n}
\mathrm{E}\Bigl(\sup_{\bolds\theta\in B_{\delta_{1}}
(\bolds\theta_{0})}
\|\nabla g_{t}(\bolds\theta)\|^{4 q_{1}}\Bigr)&=&O(1).
\end{eqnarray}

\item
For any $\delta_{2}>0$ such that
$\boldTheta_{1}-B_{\delta_{2}}(\bolds\theta_{0})$
is nonempty,
\textup{(C2)}--\textup{(C4)} hold for $\boldTheta=
\boldTheta_{1}-
B_{\delta_{2}}(\bolds\theta_{0})$,
$\mathbf{f}_{t}(\bolds\theta)= g_{t}
(\bolds\theta)- g_{t} (\bolds\theta_{0})$
and $\mathcal{F}_{t}=\mathcal{G}_{t-1}$.
In addition, there exist $0< \nu\leq1/2$
and $q_{2}>q/(2 \nu)$
such that
%
%
\begin{equation}\label{equ2.11}
\mathrm{E}\Biggl(\sup_{\bolds\theta\in\boldTheta_
{1}-B_{\delta_{2}}(\bolds\theta_{0})}
\Biggl|n^{-1}\sum_{t=1}^{n}\varepsilon_{t}\bigl(g_{t}(\bolds\theta)-
g_{t}(\bolds\theta_{0})\bigr)\Biggr|^{q_{2}}\Biggr) =O(n^{-\nu q_{2}}).
\end{equation}

\item
There exists $\bar{M}>0$
such that
%
%
\begin{eqnarray}\label{equ2.12}\qquad
P\Biggl(\sup_{\bolds\theta\in B_{\delta_{1}}
(\bolds\theta_{0})}
\lambda^{-1}_{\min}\Biggl(n^{-1}
\sum_{t=1}^{n}\nabla g_{t}(\bolds\theta) (\nabla g_{t}
(\bolds\theta))^{\mathrm T}\Biggr)
>\bar{M}
\Biggr) &=&O(n^{-q}),
\\
\label{equ2.13}
P\Biggl(\sup_{\bolds\theta\in B_{\delta_{1}}
(\bolds\theta_{0})} n^{-1}
\sum_{t=1}^{n}\|\nabla g_{t}(\bolds\theta)\|^{2}
>\bar{M}
\Biggr) &=&O(n^{-q}),
\\
\label{equ2.14}
\max_{1 \leq i, j \leq k}P\Biggl(\sup_{\bolds\theta\in
B_{\delta_{1}}(\bolds\theta_{0})} n^{-1}
\sum_{t=1}^{n}(\nabla^{2} g_{t}(\bolds\theta))_{i,j}^{2}
>\bar{M}
\Biggr) &=&O(n^{-q}).
\end{eqnarray}
\end{longlist}
Then (\ref{equ2.7}) holds.
\end{theorem}

Some comments are in order.
Conditions (i) and (iii) are needed to prove
that the $q$th moment of $\|n^{1/2}(\hat{\bolds\theta}_{n} -
\bolds\theta_{0})\|I_{A_{n}}$ is asymptotically bounded in (\ref{equ2.15}),
where $A_{n}$ is the event $\hat{\bolds\theta}_{n}$ falls into a
small ball
around $\bolds\theta_{0}$. Equations (\ref{equ2.9}) and (\ref{equ2.10})
in condition (i)
are similar to Condition 13 of~\cite{17}, but (\ref{equ2.9}) and
(\ref
{equ2.10}) require the existence of higher-order moments of
$\nabla g_{t}(\bolds\theta)$ and $\nabla^{2} g_{t}(\bolds\theta)$
to establish inequality~(\ref{equ2.26}), which plays an important role
in deriving (\ref{equ2.15}).
Equation (\ref{equ2.8}) in condition (i) can be viewed as a ``moment''
counterpart to (3.18) of~\cite{14}
and can be justified by an argument similar to (3.8) of~\cite{14},
which shows that the supremum of a Hilbert space ($H$) valued
martingale is dominated by its norm in $H$ under certain smoothness conditions.
For more details, see (\ref{equB.5}) and (\ref{equB.7}) of Appendix~\ref{appB}.
Equations (\ref{equ2.12})--(\ref{equ2.14}) in condition (iii) may seem
less relevant
to the typical assumptions made for LLN and CLT of
$\hat{\bolds\theta}_{n}$ at the first sight.
However, like (\ref{equ2.9}) and~(\ref{equ2.10}),
they are needed for the derivation of (\ref{equ2.26}).
In fact, (\ref{equ2.12}) and (\ref{equ2.13}) can be simplified into a
single assumption that
for any $\bar{m}>0$,
\begin{eqnarray*}
&&P\Biggl(\sup_{\bolds\theta\in
B_{\delta_{1}}(\bolds\theta_{0})} \Biggl\|n^{-1}
\sum_{t=1}^{n}[\nabla g_{t}(\bolds\theta)(\nabla g_{t}(\bolds\theta
))^{\mathrm T}-
\mathrm{E} \{\nabla g_{t}(\bolds\theta)(\nabla g_{t}(\bolds\theta
))^{\mathrm T}\}]\Biggr\|
>\bar{m}
\Biggr) \\
&&\qquad=O(n^{-q}),
\end{eqnarray*}
where $\|D\|^{2}=\sup_{\|\mathbf{x}\|=1} \mathbf{x}^{\mathrm
T}D^{\mathrm T}D\mathbf{x}$
for the matrix $D$.
However, we do not want to complicate the proof of Theorem~\ref{theo2.2}
by using this assumption.
When $g_{t}(\bolds\theta)$
is a linear process with coefficient functions satisfying certain
smoothness conditions,
(\ref{equ2.12})--(\ref{equ2.14}) can be justified based on
a \textit{uniform} version of the first moment bound
theorem of Findley and Wei~\cite{6}. Further details can be found in
(\ref{equB.6}) and (\ref{equB.9})--(\ref{equB.11}) of
Appendix~\ref{appB}. In contrast to conditions (i) and (iii), condition (ii) is
required to prove that
the $q$th moment of $\|n^{1/2}(\hat{\bolds\theta}_{n} -
\bolds\theta_{0})\|I_{B_{n}}$ is asymptotically bounded in (\ref{equ2.27}),
where $B_{n}$ denotes the event
$\hat{\bolds\theta}_{n}$ falls \textit{outside} a small ball
around~$\bolds\theta_{0}$. Finally, (C2) in condition (ii)
provides an identifiability condition for model (\ref{equ2.4}), while
(\ref{equ2.11}) is a moment counterpart to (3.14) of~\cite{14}
and can be analogously justified as (\ref{equ2.8}).
\begin{pf*}{Proof of Theorem~\ref{theo2.2}}
Let $0< \delta^{*}_{1} < \min\{\delta_{1},
3^{-1}k^{-1}\bar{M}^{-2}\}$ and $A_{n}=\{\hat{\bolds\theta}_{n} \in
B_{\delta_{1}^{*}} (\bolds\theta_{0})\}$. We first show that
%
%
\begin{equation}\label{equ2.15}
\mathrm{E}\bigl(\|n^{1/2}(\hat{\bolds\theta}_{n} -
\bolds\theta_{0})\|^{q}I_{A_{n}}\bigr) = O(1).
\end{equation}

By the mean value theorem for vector-valued functions,
on the set $A_{n}$,
%
%
\begin{equation}\label{equ2.16}
\mathbf{0}= \nabla S_{n}(\hat{\bolds\theta}_{n})
=\nabla S_{n}(\bolds\theta_{0})+
\biggl\{\!\int_{0}^{1}\! \nabla^{2} S_{n}\bigl(\bolds\theta_{0}+
r(\hat{\bolds\theta}_{n}-\bolds\theta_{0})\bigr) \,dr
\!\biggr\}
(\hat{\bolds\theta}_{n}-\bolds\theta_{0}),\hspace*{-35pt}
\end{equation}
where $S_{n}(\cdot)$ is defined in (\ref{equ2.6})
and
the integral of a matrix is to be understood component-wise.
In view of (\ref{equ2.16}) and the identities that
$\nabla S_{n}(\bolds\theta)
=-2 \sum_{t=1}^{n}(y_{t}-g_{t}(\bolds\theta)) \nabla g_{t}
(\bolds\theta)$ and
$\nabla^{2} S_{n}(\bolds\theta)=
2 \sum_{t=1}^{n} \nabla g_{t}(\bolds\theta) (\nabla g_{t}
(\bolds\theta))^{\mathrm T}
-\break 2 \sum_{t=1}^{n}(y_{t}-g_{t}(\bolds\theta))\times \nabla^{2} g_{t}
(\bolds\theta)$,
one has
%
%
\begin{equation}\label{equ2.17}\quad
\sum_{t=1}^{n} \varepsilon_{t} \nabla g_{t}(\bolds\theta_{0})
=\bigl(L(\hat{\bolds\theta}_{n}, \bolds\theta_{0})
-Q(\hat{\bolds\theta}_{n}, \bolds\theta_{0})\bigr)
(\hat{\bolds\theta}_{n} - \bolds\theta_{0})
\qquad\mbox{on } A_{n},
\end{equation}
where
$L(\hat{\bolds\theta}_{n}, \bolds\theta_{0})=
\int_{0}^{1}
\sum_{t=1}^{n} \nabla g_{t}(\bolds\theta_{0}+
r(\hat{\bolds\theta}_{n} - \bolds\theta_{0}))
(\nabla g_{t}(\bolds\theta_{0}+
r(\hat{\bolds\theta}_{n} - \bolds\theta_{0})))^
{\mathrm T} \,dr$ and
$Q(\hat{\bolds\theta}_{n}, \bolds\theta_{0})=
\int_{0}^{1}
\sum_{t=1}^{n}
\{y_{t}-g_{t}(\bolds\theta_{0}+
r(\hat{\bolds\theta}_{n} - \bolds\theta_{0}))\}
\nabla^{2} g_{t}(\bolds\theta_{0}+
r(\hat{\bolds\theta}_{n} - \bolds\theta_{0})) \,dr$.
A direct algebraic manipulation leads to
%
%
\begin{equation}\label{equ2.18}
\lambda_{\min}(L(\hat{\bolds\theta}_{n},
\bolds\theta_{0}))
\geq\inf_{\bolds\theta\in B_{\delta_{1}}
(\bolds\theta_{0})}
\lambda_{\min}\Biggl(
\sum_{t=1}^{n} \nabla g_{t}(\bolds\theta) (\nabla g_{t}
(\bolds\theta))^{\mathrm T}\Biggr) \qquad\mbox{on } A_{n},\hspace*{-35pt}
\end{equation}
which, together with the continuity of $\nabla g_{t}(\bolds\theta)$ on
$B_{\delta_{1}}(\bolds\theta_{0})$,
condition (i) and Theorem~\ref{theo2.1}, yields that for any $s \geq1$,
%
%
\begin{eqnarray}\label{equ2.19}\quad
&&\mathrm{E}(\lambda^{-s}_{\min}(n^{-1}L(\hat{\bolds\theta}_{n},
\bolds\theta_{0}))I_{A_{n}}) \nonumber\\[-8pt]\\[-8pt]
&&\qquad\leq\mathrm{E}\Biggl( \sup_{\bolds\theta\in B_{\delta_{1}}
(\bolds\theta_{0})}
\lambda^{-s}_{\min}\Biggl(n^{-1}\sum_{t=1}^{n} \nabla g_{t}(\bolds\theta)
(\nabla g_{t}(\bolds\theta))^{\mathrm T}\Biggr) \Biggr)
= O(1).\nonumber
\end{eqnarray}
With the help of (\ref{equ2.19}), we can assume without loss of
generality that
$L^{-1}(\hat{\bolds\theta}_{n}, \bolds\theta_{0})$
exists on $A_{n}$, and hence by (\ref{equ2.17}),
%
%
\begin{eqnarray}\label{equ2.20}
&&\|n^{1/2}(\hat{\bolds\theta}_{n} -
\bolds\theta_{0})\|I_{A_{n}} \nonumber\\
&&\qquad\leq
\|nL^{-1}(\hat{\bolds\theta}_{n}, \bolds\theta_{0})\|
\Biggl\|n^{-1/2}\sum_{t=1}^{n}\varepsilon_{t} \nabla g_{t}
(\bolds\theta_{0})\Biggr\|I_{A_{n}} \nonumber\\
&&\qquad\quad{}+ \|nL^{-1}(\hat{\bolds\theta}_{n},
\bolds\theta_{0})\|\Biggl\|
\int_{0}^{1} n^{-1/2}
\sum_{t=1}^{n}\varepsilon_{t}
\nabla^{2} g_{t}\bigl(\bolds\theta_{0}+
r(\hat{\bolds\theta}_{n} - \bolds\theta_{0})\bigr) \,dr
\Biggr\|\nonumber\\
&&\hspace*{44pt}{}\times \|\hat{\bolds\theta}_{n}- \bolds\theta_{0}\|I_{A_{n}} \\
&&\qquad\quad{}+
\|nL^{-1}(\hat{\bolds\theta}_{n}, \bolds\theta_{0})\|\Biggl\|
\int_{0}^{1} n^{-1}
\sum_{t=1}^{n} r (\bolds\theta_{0}-
\hat{\bolds\theta}_{n})^{\mathrm T}
\nabla g_{t}(\bolds\theta^{*}_{t, r})\nonumber\\
&&\hspace*{169.3pt}{}\times\nabla^{2} g_{t}\bigl(\bolds\theta_{0}+r(\hat{\bolds\theta}
_{n} - \bolds\theta_{0})\bigr) \,dr \Biggr\| \nonumber\\
&&\hspace*{44pt}{}\times\|n^{1/2}(\hat{\bolds\theta}_{n}-
\bolds\theta_{0})\|I_{A_{n}},\nonumber
\end{eqnarray}
where $\bolds\theta^{*}_{t, r}$ satisfies
$\|\bolds\theta^{*}_{t, r}-\bolds\theta_{0}\| \leq
r\|\hat{\bolds\theta}_{n} - \bolds\theta_{0}\|$.
By the Cauchy--Schwarz inequality and Jensen's inequality,
it follows that
%
%
\begin{eqnarray}\label{equ2.21}
&&\Biggl\|
\int_{0}^{1} n^{-1/2}
\sum_{t=1}^{n}\varepsilon_{t}
\nabla^{2} g_{t}\bigl(\bolds\theta_{0}+
r(\hat{\bolds\theta}_{n} -
\bolds\theta_{0})\bigr) \,dr \Biggr\| \nonumber\\[-8pt]\\[-8pt]
&&\qquad\leq
k \max_{1 \leq i, j \leq k} \sup_{\bolds\theta\in
B_{\delta_{1}}(\bolds\theta_{0})}
\Biggl|n^{-1/2}\sum_{t=1}^{n}\varepsilon_{t}
(\nabla^{2} g_{t}(\bolds\theta))_{i, j}\Biggr| := k W_{n}\nonumber
\end{eqnarray}
and
%
%
\begin{eqnarray}\label{equ2.22}
&& \Biggl\|
\int_{0}^{1} n^{-1}
\sum_{t=1}^{n} r (\bolds\theta_{0}-
\hat{\bolds\theta}_{n})^{\mathrm T}
\nabla g_{t}(\bolds\theta^{*}_{t, r})
\nabla^{2} g_{t}\bigl(\bolds\theta_{0}+r(\hat{\bolds\theta}
_{n} - \bolds\theta_{0})\bigr) \,dr \Biggr\| \nonumber\\
&&\hspace*{2.8pt}\qquad\leq k \|\hat{\bolds\theta}_{n}-\bolds\theta_{0}\|
\Biggl(\sup_{\bolds\theta\in B_{\delta_{1}}
(\bolds\theta_{0})}
n^{-1}\sum_{t=1}^{n}\|\nabla g_{t}(\bolds\theta)\|^{2}
\Biggr)^{1/2} \nonumber\\[-8pt]\\[-8pt]
&&\qquad\quad\hspace*{0pt}{}\times
\Biggl\{\max_{1 \leq i, j \leq k}\sup_{\bolds\theta\in
B_{\delta_{1}}(\bolds\theta_{0})}
n^{-1}\sum_{t=1}^{n}(\nabla^{2} g_{t}(\bolds\theta))_{i, j}^{2}
\Biggr\}^{1/2}\nonumber\\
&&\qquad:= k\|\hat{\bolds\theta}_{n}-\bolds\theta_{0}\|
R^{1/2}_{1n}R^{1/2}_{2n}.\nonumber
\end{eqnarray}
Denoting $\sup_{\bolds\theta\in B_{\delta_{1}}
(\bolds\theta_{0})}
\lambda^{-1}_{\min}(n^{-1} \sum_{t=1}^{n} \nabla
g_{t}(\bolds\theta)
(\nabla g_{t}(\bolds\theta))^{\mathrm T})$
by $R_{n}$ and combining (\ref{equ2.18}) and (\ref{equ2.20})--(\ref
{equ2.22}), we obtain
%
%
\begin{eqnarray}\label{equ2.23}\qquad
&&\|n^{1/2}(\hat{\bolds\theta}_{n} -
\bolds\theta_{0})\|^{q}I_{A_{n}} \nonumber\\
&&\hspace*{2.8pt}\qquad\leq 3^{q}\Biggl\{R^{q}_{n}
\Biggl\|n^{-1/2}\sum_{t=1}^{n}\varepsilon_{t} \nabla g_{t}
(\bolds\theta_{0})\Biggr\|^{q} \nonumber\\[-8pt]\\[-8pt]
&&\qquad\quad\hspace*{18.4pt}{} + \delta^{*^{q}}_{1}k^{q}R^{q}_{n}W^{q}_{n} +
\delta^{*^{q}}_{1}k^{q}R^{q}_{n} R^{q/2}_{1n}R^{q/2}_{2n}
\|n^{1/2}(\hat{\bolds\theta}_{n}- \bolds\theta_{0})\|
^{q}I_{A_{n}}\Biggr\} \nonumber\\
&&\qquad:= 3^{q}\{\mathrm{(I)}+\mathrm{(II)}+\mathrm{(III)}\}.\nonumber
\end{eqnarray}
Applying (\ref{equ2.8}), (\ref{equ2.10}), (\ref{equ2.19}), $\sup
_{t}\mathrm{E}(|\varepsilon_{t}|^{\gamma}|\mathcal{G}_{t-1})<
C_{5}$ a.s., H\"{o}lder's inequality and Lemma 2 of Wei~\cite{22},
it can be shown that for $n$ large and some positive constants
$C^{*}_{1}$ and $C^{*}_{2}$,
%
%
\begin{equation}\label{equ2.24}
\mathrm{E}(\mathrm{(I)}) \leq C^{*}_{1}
\end{equation}
and
%
%
\begin{equation}\label{equ2.25}
\mathrm{E}(\mathrm{(II)}) \leq C^{*}_{2};
\end{equation}
see Appendix A of Chan and Ing~\cite{4} for more details.
In addition, by making use of
(\ref{equ2.9}), (\ref{equ2.10}) and (\ref{equ2.12})--(\ref{equ2.14}),
we show in Appendix~\ref{appA} that for $n$ large,
%
%
\begin{equation}\label{equ2.26}
\mathrm{E}(\mathrm{(III)}) \leq C^{*}_{3}+C^{*}_{4}
\mathrm{E}\bigl(\|n^{1/2}(\hat{\bolds\theta}_{n} -
\bolds\theta_{0})\|^{q}I_{A_{n}}\bigr),
\end{equation}
where $C^{*}_{3}$ and $C^{*}_{4}$ are some positive constants
with $C^{*}_{4}$
satisfying $0<C^{*}_{4}<3^{-q}$.
Consequently,
the desired conclusion (\ref{equ2.15}) follows from (\ref
{equ2.23})--(\ref{equ2.26}).

Letting $B_{n}=\{\hat{\bolds\theta}_{n} \in
\tilde{\boldTheta}_{1}=
\boldTheta_{1}-B_{\delta^{*}_{1}}
(\bolds\theta_{0})\}$,
the rest of the proof aims to show that
%
%
\begin{equation}\label{equ2.27}
\mathrm{E}\bigl(\|n^{1/2}(\hat{\bolds\theta}_{n} -
\bolds\theta_{0})\|^{q}I_{B_{n}}\bigr) = O(1),
\end{equation}
which, together with (\ref{equ2.15}), yields the desired conclusion
(\ref{equ2.7}).

Since $\|n^{1/2}(\hat{\bolds\theta}_{n} -
\bolds\theta_{0})\|^{q} \leq C^{*}_{5}n^{q/2}$
for some $C^{*}_{5}>0$, (\ref{equ2.27}) follows immediately once we can
show that
%
%
\begin{equation}\label{equ2.28}
P(B_{n}) = O(n^{-q/2}).
\end{equation}
By the continuity of $g_{t}(\cdot)$ on $\boldTheta_{1}$,
condition (ii) and Theorem~\ref{theo2.1}, one has for any $s \geq1$,
%
%
\begin{equation}\label{equ2.29}
\mathrm{E}\Biggl\{
\Biggl[\inf_{\bolds\theta\in\tilde{\boldTheta}_{1}}
n^{-1}\sum_{t=1}^{n}\bigl(g_{t}(\bolds\theta)-g_{t}
(\bolds\theta_{0})\bigr)^{2}\Biggr]^{-s}\Biggr\}=O(1).
\end{equation}
In addition, it is straightforward to see that
%
%
\begin{eqnarray}\label{equ2.30}
&&
B_{n} \subseteq\Biggl\{2 \sup_{\bolds\theta\in
\tilde{\boldTheta}_{1}}
\Biggl|n^{-1}\sum_{t=1}^{n}\varepsilon_{t}\bigl(g_{t}(\bolds\theta)-g_{t}
(\bolds\theta_{0})\bigr)\Biggr|\nonumber\\[-8pt]\\[-8pt]
&&\qquad\hspace*{15.8pt}\geq\inf_{\bolds\theta\in\tilde{\boldTheta}_{1}}
n^{-1}\sum_{t=1}^{n}\bigl(g_{t}(\bolds\theta)-g_{t}
(\bolds\theta_{0})\bigr)^{2}\Biggr\}.\nonumber
\end{eqnarray}
Since $q_{2}>q/(2 \nu)$,
there exists $\eta_{1}>0$ such that $q_{2}=q(1+\eta_{1})/(2 \nu)$.
By (\ref{equ2.29}), (\ref{equ2.30}), (\ref{equ2.11}), Chebyshev's
inequality and H\"{o}lder's inequality,
there exists $C^{*}_{6}>0$ such that for all large $n$,
\begin{eqnarray*}
P(B_{n}) &\leq& C^{*}_{6}
\Biggl\{\mathrm{E} \Biggl(\inf_{\bolds\theta\in
\tilde{\boldTheta}_{1}}
n^{-1}\sum_{t=1}^{n}\bigl(g_{t}(\bolds\theta)-g_{t}
(\bolds\theta_{0})\bigr)^{2}\Biggr)^{-q_{2}/\eta_{1}}
\Biggr\}^{\eta_{1}/(1+\eta_{1})} \\
&&\hspace*{0pt}{}\times\Biggl\{\mathrm{E} \Biggl(\sup_{\bolds\theta\in
\tilde{\boldTheta}_{1}}
\Biggl|n^{-1}\sum_{t=1}^{n}
\varepsilon_{t}\bigl(g_{t}(\bolds\theta)-g_{t}
(\bolds\theta_{0})\bigr)\Biggr|^{q_{2}}\Biggr)\Biggr\}^{1/(1+\eta_{1})}= O(n^{-q/2}).
\end{eqnarray*}
Consequently, (\ref{equ2.28}) is established and the theorem is proved.
\end{pf*}

As mentioned in the \hyperref[intro]{Introduction},
(\ref{equ2.7}) can be used to examine the
asymptotic properties of MSPE
of $g_{n+1}(\hat{\bolds\theta}_{n})$, ${\mathrm
E}(y_{n+1}-g_{n+1}(\hat
{\bolds\theta}_{n}))^{2}$,
which is also known as the \textit{final prediction error} (FPE) for AR
models; see Akaike~\cite{1}.
To see this, note first that under certain mild conditions such as
(2.2) and (2.3) of~\cite{14},
$\hat{\bolds\theta}_{n} \to\bolds\theta_{0}$ a.s. If one can
further show that
%
%
\begin{equation}\label{equ2.31}
n^{1/2} (\nabla g_{n+1}(\bolds\theta_{0}))^{\mathrm T}(\hat{\bolds
\theta
}_{n}-\bolds\theta_{0})
\Rightarrow\mathbf{H}
\end{equation}
and
%
\begin{equation}
\label{equ2.32}
n\{(\nabla g_{n+1}(\bolds\theta_{0}))^{\mathrm T}(\hat{\bolds
\theta
}_{n}-\bolds\theta_{0})\}^{2} \mbox{ is uniformly
integrable},
\end{equation}
where $\Rightarrow$ denotes convergence in distribution and $\mathbf
{H}$ is a random variable
with ${\mathrm E}(\mathbf{H}^{2})< \infty$, then
%
%
\begin{equation}\label{equ2.33}
\lim_{n \to\infty} n {\mathrm E}
\{(\nabla g_{n+1}(\bolds\theta_{0}))^{\mathrm T}(\hat{\bolds\theta
}_{n}-\bolds\theta_{0})\}^{2} =
{\mathrm E} (\mathbf{H}^{2}).
\end{equation}
Once (\ref{equ2.33}) is established, it can be linked to
${\mathrm E}(y_{n+1}-g_{n+1}(\hat{\bolds\theta}_{n}))^{2}$
by means of Taylor's expansion as follows. Note that
%
%
\begin{eqnarray}\label{equ2.34}
&& n\bigl\{{\mathrm E}\bigl(y_{n+1}-g_{n+1}(\hat{\bolds\theta}_{n})\bigr)^{2}-
\mathrm
{E}(\varepsilon^{2}_{n+1})\bigr\} \nonumber\\[-8pt]\\[-8pt]
&&\qquad=
n{\mathrm E}\{(\nabla g_{n+1}(\bolds\theta_{0}))^{\mathrm T}(\hat
{\bolds\theta}_{n}-\bolds\theta_{0})\}^{2}+
{\mathrm E}(\tilde{R}_n) \to{\mathrm E} (\mathbf{H}^{2}),
\nonumber
\end{eqnarray}
provided the remainder term $\tilde{R}_n$ satisfies E($\tilde
{R}_n$)$=o(1)$. While (\ref{equ2.31}) can be
established by means of asymptotic distribution results (see Section~\ref{sec3}),
(\ref{equ2.7}) serves as an important device in
establishing (\ref{equ2.32}) and E($\tilde{R}_n$)$=o(1)$. If one
further assumes that
$\mathrm{E}(\varepsilon^{2}_{t})=\sigma^{2}>0$ for all $t>0$, then
(\ref
{equ2.34}) provides an asymptotic expression for
${\mathrm E}(y_{n+1}-g_{n+1}(\hat{\bolds\theta}_{n}))^{2}$ as
%
%
\begin{equation}\label{equ2.35}
{\mathrm E}\bigl(y_{n+1}-g_{n+1}(\hat{\bolds\theta}_{n})\bigr)^{2}=\sigma
^{2}+\frac
{{\mathrm E} (\mathbf{H}^{2})}{n}+o(n^{-1}).
\end{equation}
Although the second term in (\ref{equ2.35}) is asymptotically
negligible compared to $\sigma^{2}$,
${\mathrm E} (\mathbf{H}^{2})$ becomes a key quantity. Utilizing (\ref
{equ2.7}), one can
make use of the asymptotic expression in (\ref{equ2.35}), in particular
E($\mathbf{H}^{2}$), to construct optimal model selection
criteria; see, for example, Akaike~\cite{1}, Wei~\cite{23} and Findley
and Wei~\cite{7}.
See, also, Section~\ref{sec3} for further discussions.

\section{Applications to ARMA models}\label{sec3}

Let $y_{1}, \ldots, y_{n}$ be
generated from the stochastic regression model,
%
%
\begin{equation}\label{equ3.1}
y_{t}=g_{t}(\boldeta_{0})+ \varepsilon_{t},\qquad t=1, \ldots, n,
\end{equation}
where $\boldeta_{0}=(\alpha_{0,1},
\ldots,\alpha_{0, p_{1}}, \beta_{0,1}, \ldots, \beta
_{0,p_{2}})^{\mathrm T}$
is an unknown coefficient vector and
$g_{t}(\boldeta_{0})$
has the ARMA representation
%
%
\begin{equation}\label{equ3.2}\quad
g_{t}(\boldeta_{0})=\alpha_{0,1}y_{t-1}+ \cdots+ \alpha_{0,
p_{1}}y_{t-p_{1}}-\beta_{0,1}\varepsilon_{t-1} - \cdots-
\beta_{0,p_{2}} \varepsilon_{t-p_{2}}
\end{equation}
with the initial conditions $y_{t}=\varepsilon_{t}=0$ for all $t\leq
0$.
Define
\[
\hat{\boldeta}_{n}=\mathop{\arg\min}_{\boldeta\in\Pi} \sum_{t=1}^{n}
\bigl(y_{t}-g_{t}(\boldeta)\bigr)^{2},
\]
where $\Pi\subset R^{p_{1}+p_{2}}$
is a compact set that includes $\boldeta_{0}$
as an interior point and
whose
elements
$\boldeta=(\alpha_{1}, \ldots, \alpha_{p_{1}}, \beta_{1}, \ldots,
\beta
_{p_{2}})^{\mathrm T}$
satisfy the following properties:
%
%
\begin{eqnarray}\label{equ3.3}
&& A_{1, \boldeta}(z) = 1-\sum_{j=1}^{p_{1}}\alpha_{j}z^{j}
\neq0,\nonumber\\[-8pt]\\[-8pt]
&&A_{2, \boldeta}(z) = 1-\sum_{j=1}^{p_{2}}\beta_{j}z^{j} \neq0 \qquad\mbox
{for all } |z| \leq1;\nonumber\\
\label{equ3.4}
&& A_{1, \boldeta}(z) \mbox{ and }
A_{2, \boldeta}(z) \mbox{ have no common zeros}; \\
\label{equ3.5}
&& |\alpha_{p_{1}}|+|\beta_{p_{2}}|>0.
\end{eqnarray}
In this section, we apply the results obtained in Section~\ref{sec2}
to show that
%
%
\begin{equation}\label{equ3.6}
\mathrm{E} \|n^{1/2}(\hat{\boldeta}_{n}-\boldeta_{0})\|^{q}
=O(1),\qquad
q\geq1.
\end{equation}
Applications of (\ref{equ3.6}) to the investigation of the MSPE
of $g_{n+1}(\hat{\boldeta}_{n})$,
$\mathrm{E} (y_{n+1}-g_{n+1}(\hat{\boldeta}_{n}))^{2}$,
are also given.
It should be mentioned that our initial conditions, $y_{t}=\varepsilon
_{t}=0$ for all\vadjust{\goodbreak} $t\leq
0$, are made for simplicity of the argument only and
all results in this section can be straightforwardly extended to the
case where
$(y_{t}, \varepsilon_{t})$ obey the same assumptions for $t \leq0$ as
for $t > 0$.

Let $\boldeta\in\Pi$.
Define
$\varepsilon_{t}(\boldeta) = 0$ for $t \leq0$ and define
$\varepsilon_{t} (\boldeta)$ recursively for $t \geq1$ by
%
%
\begin{eqnarray}\label{equ3.7}
\varepsilon_{t} (\boldeta)&=&y_{t}-g_{t}(\boldeta)
\nonumber\\
&=& y_{t}-\alpha_{1}y_{t-1}-\cdots-\alpha_{p}y_{t-p_{1}}\\
&&{}+\beta_{1}
\varepsilon_{t-1}(\boldeta) + \cdots+
\beta_{p_{2}}\varepsilon_{t-p_{2}}(\boldeta),\nonumber
\end{eqnarray}
noting that $\varepsilon_{t} (\boldeta_{0})=\varepsilon_{t}$.
As observed in (\ref{equ2.19}) and (\ref{equ2.29}) of Section~\ref{sec2},
to obtain (\ref{equ3.6}),
it is crucial to verify that for some $\delta_{1} >0$ with
$B_{\delta_{1}}(\boldeta_{0}) \subset\Pi$
and any $s\geq1$,
%
%
\begin{equation}\label{equ3.8}
\mathrm{E}\Biggl\{\sup_{\boldeta\in B_{\delta_{1}}(\boldeta_{0})}
\lambda^{-s}_{\min}\Biggl[n^{-1}\sum_{t=1}^{n} \nabla\varepsilon
_{t}(\boldeta)
(\nabla\varepsilon_{t}(\boldeta))^{\mathrm T}\Biggr]\Biggr\}=O(1);
\end{equation}
and for any $\delta_{2}>0$ with
$\tilde{\Pi}=\Pi-B_{\delta_{2}}(\boldeta_{0}) \neq\varnothing$
and any $s\geq1$,
%
%
\begin{equation}\label{equ3.9}
\mathrm{E}\Biggl\{
\sup_{\boldeta\in\tilde{\Pi}}
\Biggl[n^{-1}\sum_{t=1}^{n}\bigl(\varepsilon_{t}(\boldeta)-
\varepsilon_{t}(\boldeta_{0})\bigr)^{2}\Biggr]^{-s}
\Biggr\}=O(1).
\end{equation}
Denote the $i$th component of $\nabla\varepsilon_{t}(\boldeta)$
by $(\nabla\varepsilon_{t}(\boldeta))_{i}$.
Straightforward calculations yield that for $1 \leq i \leq p_{1}$
and $1 \leq j \leq p_{2}$,
%
%
\begin{eqnarray}\label{equ3.10}
(\nabla\varepsilon_{t} (\boldeta))_{i} &=&
- y_{t-i} +
\sum_{s=1}^{p_{2}} \beta_{s}
(\nabla\varepsilon_{t-s} (\boldeta))_{i},\\
\label{equ3.11}
(\nabla\varepsilon_{t} (\boldeta))_{p_{1}+j} &=&
\varepsilon_{t-j}(\boldeta) +
\sum_{s=1}^{p_{2}} \beta_{s}
(\nabla\varepsilon_{t-s} (\boldeta))_{p_{1}+j}.
\end{eqnarray}
For $j<0$, let $c_{j}^{(1)}(\boldeta)=c_{j}^{(2)}(\boldeta)=0$ and
for $j \geq0$, let $c_{j}^{(1)}(\boldeta)$ and $c_{j}^{(2)}(\boldeta)$
satisfy
%
%
\begin{equation}\label{equ3.12}
\sum_{j=0}^{\infty}c_{j}^{(1)}(\boldeta)z^{j} =
\frac{-A_{2, \boldeta_{0}}(z)}
{A_{2, \boldeta}(z)A_{1, \boldeta_{0}}(z)},\qquad
\sum_{j=0}^{\infty}c_{j}^{(2)}(\boldeta)z^{j} =
\frac{A_{1, \boldeta}(z)A_{2, \boldeta_{0}}(z)}
{A^{2}_{2, \boldeta}(z)A_{1, \boldeta_{0}}(z)}.\hspace*{-35pt}
\end{equation}
In view of (\ref{equ3.3})--(\ref{equ3.5}) and the compactness of $\Pi$,
there exist
positive constants $K_{1}$ and $K_{2}$ such that for all $j\geq0$
and $i=1, 2$,
%
%
\begin{equation}\label{equ3.13}
\sup_{\boldeta\in\Pi}\bigl|c^{(i)}_{j}(\boldeta)\bigr| \leq K_{1}
\operatorname{exp}(-K_{2}j).
\end{equation}
Define
$b^{(l)}_{j}(\boldeta)=c_{j-l}^{(1)}(\boldeta), 1 \leq l \leq p_{1}$, and
$b^{(p_{1}+l)}_{j}(\boldeta)=c_{j-l}^{(2)}(\boldeta), 1 \leq l \leq p_{2}$.
Then it follows from (\ref{equ3.10})--(\ref{equ3.13}) that
%
%
\begin{equation}\label{equ3.14}
\nabla\varepsilon_{t}(\boldeta) =
\Biggl(\sum_{j=1}^{t-1}b^{(1)}_{j}(\boldeta)\varepsilon_{t-j}, \ldots,
\sum_{j=1}^{t-1}b^{(p_{1}+p_{2})}_{j}(\boldeta)\varepsilon
_{t-j}\Biggr)^{\mathrm T}
\end{equation}
and
\begin{equation}
\label{equ3.15}
\max_{1 \leq l \leq p_{1}+p_{2}}\sup_{\boldeta\in\Pi
}\bigl|b^{(l)}_{j}(\boldeta)\bigr| \leq K^{\prime}_{1} \operatorname{exp}(-K_{2}j)\qquad \mbox
{for some }
K^{\prime}_{1}>0.
\end{equation}
Moreover, one has
%
%
\begin{equation}\label{equ3.16}
\varepsilon_{t}(\boldeta)-\varepsilon_{t}(\boldeta_{0})=
\sum_{i=1}^{t-1}b_{i}(\boldeta)\varepsilon_{t-i},
\end{equation}
where $b_{j}(\boldeta), j \geq1$, satisfy
$1+\sum_{j=1}^{\infty}b_{j}(\boldeta)z^{j}=
A_{1, \boldeta}(z)A_{2, \boldeta_{0}}(z)/
(A_{2, \boldeta}(z)\times A_{1, \boldeta_{0}}(z))$
and
%
%
\begin{equation}\label{equ3.17}
\sup_{\boldeta\in\Pi}|b_{j}(\boldeta)|\leq K_{3} \operatorname{exp}(-K_{4}j)
\end{equation}
for some positive constants $K_{3}$ and $K_{4}$.
The next
theorem provides sufficient conditions under which
%
%
\begin{equation}\label{equ3.18}
\mathrm{E}\Biggl\{\sup_{\boldeta\in\Pi}
\lambda^{-s}_{\min}\Biggl[n^{-1}\sum_{t=1}^{n} \nabla\varepsilon
_{t}(\boldeta)
(\nabla\varepsilon_{t}(\boldeta))^{\mathrm T}\Biggr]\Biggr\}=O(1)\qquad \mbox{for any }
s\geq1.\hspace*{-35pt}
\end{equation}
This result leads immediately to (\ref{equ3.8}).
\begin{theorem}\label{theo3.1}
Assume model (\ref{equ3.1}), with $g_{t}(\cdot)$ defined in
(\ref{equ3.2})
and $\varepsilon_{t}$ being
independent random variables satisfying $\mathrm{E}(\varepsilon_{t})=0$
and $\mathrm{E}(\varepsilon^{2}_{t})=\sigma^{2}$
for all $t \geq1$.
Moreover, assume that there exist positive constants $\alpha_{1}, \xi$
and $M_{1}$ such that for any
$0 < s_{2} - s_{1} \leq\xi$,
%
%
\begin{equation}\label{equ3.19}
{\sup_{1 \leq m \leq t < \infty, \|\mathbf{v}\|=1}}
|
F_{t, m, \mathbf{v}}(s_{2}) -F_{t, m, \mathbf{v}}(s_{1}) |
\leq M_{1}(s_{2} - s_{1})^{\alpha_{1}},
\end{equation}
where $\mathbf{v} \in R^{m}$ and $F_{t, m,\mathbf{v}} (\cdot)$
denotes the distribution function of $\mathbf{v}^{\mathrm
T}(\varepsilon_{t}, \ldots,\break \varepsilon_{t+1-m})^{\mathrm T}$. Then,
\textup{(C1)}--\textup{(C4)} hold for $\boldTheta=\Pi$,
$\mathbf{f}_{t}(\bolds\theta)=\nabla\varepsilon_{t}(\boldeta)$ and
$\mathcal{F}_{t}=\sigma\{\varepsilon_{t-1},\break  \varepsilon_{t-2},
\ldots\}
$, the $\sigma$-field generated by
$\varepsilon_{t-1},\varepsilon_{t-2}, \ldots.$ Hence, by Theorem~\ref{theo2.1},
(\ref{equ3.18}) follows.
\end{theorem}
\begin{pf}
According to (\ref{equ3.3}) and (\ref{equ3.12}), it is easy to see that
$\nabla\varepsilon_{t}(\boldeta)$
is continuous on $\Pi$, and hence (C1) follows.
Define $\bolds\Lambda= \{\mathbf{a}\dvtx \mathbf{a} \in R^{\bar
{p}}, \|\mathbf{a}\|=1\}$,
where $\bar{p}=p_{1}+p_{2}$.
To show (C2), note first that by (\ref{equ3.3})--(\ref{equ3.5}),
one has
for any $\bolds\lambda\in\bolds\Lambda$ and
$\boldeta\in\Pi$,
there exists $\delta_{2}=\delta_{2}(\bolds\lambda, \boldeta
)>0$ such that
for all large $t$,
%
%
\begin{equation}\label{equ3.20}
E(\bolds\lambda^{\mathrm T} \nabla\varepsilon_{t}(\boldeta
))^{2} > \delta_{2}.
\end{equation}
In addition, it follows from (\ref{equ3.14}) and (\ref{equ3.15}) that
%
%
\begin{equation}\label{equ3.21}
\mathrm{E}(\bolds\lambda^{\mathrm T}\nabla\varepsilon
_{t}(\boldeta))^{2}
\mbox{ converges to } l(\bolds\lambda, \boldeta) \mbox{
uniformly on } \bolds\Lambda\times
\Pi,
\end{equation}
where $l(\bolds\lambda, \boldeta)$ is some nonnegative
function on $\bolds\Lambda\times\Pi$.
Moreover, since\break
$\mathrm{E}(\bolds\lambda^{\mathrm T}\nabla\varepsilon
_{t}(\boldeta))^{2}$ is continuous on
$\bolds\Lambda\times\Pi$, uniform convergence implies
that $l(\bolds\lambda, \boldeta)$ is also continuous on
$\bolds\Lambda\times\Pi$. By (\ref{equ3.20}) and the
compactness of $\bolds\Lambda\times\Pi$,
$
\inf_{\bolds\lambda\in\bolds\Lambda, \boldeta
\in\Pi}
l(\bolds\lambda, \boldeta) > 0$.
This, together with (\ref{equ3.21}), yields that
there is a positive number $\epsilon$ and a positive integer $L$ such that
for all $t > L$,
%
%
\begin{equation}\label{equ3.22}
\inf_{\bolds\lambda\in\bolds\Lambda, \boldeta
\in\Pi}
\mathrm{E}(\bolds\lambda^{\mathrm T}\nabla\varepsilon
_{t}(\boldeta))^{2}> \epsilon> 0.
\end{equation}
For $t > l_{1} \geq1$,
define $\nabla\varepsilon_{t, l_{1}}
(\boldeta) = (\sum_{i=1}^{l_{1}}
b_{i}^{(1)}(\boldeta)\varepsilon_{t-i}, \ldots, \sum_{i=1}^{l_{1}}
b_{i}^{(\bar{p})}(\boldeta)\varepsilon_{t-i})^{\mathrm T}$.
According to (\ref{equ3.14}) and (\ref{equ3.15}), there exists a
positive integer
$L_{1}(\epsilon)$ such that for all $t>l_{1}
\geq L_{1}(\epsilon)$,
%
%
\begin{equation}\label{equ3.23}
\sup_{\bolds\lambda
\in\bolds\Lambda, \boldeta\in\Pi}
|\mathrm{E}(\bolds\lambda^{\mathrm T}\nabla\varepsilon
_{t}(\boldeta))^{2} -
\mathrm{E}(\bolds\lambda^{\mathrm T}\nabla\varepsilon_{t,
l_{1}}(\boldeta))^{2}| < \epsilon/2.
\end{equation}
From (\ref{equ3.22}) and (\ref{equ3.23}), it follows that for all $t >
d_{1}=\max\{L, L_{1}(\epsilon)\}$,
%
%
\begin{equation}\label{equ3.24}
\inf_{\bolds\lambda\in\bolds\Lambda, \boldeta
\in\Pi}
\mathrm{E}(\bolds\lambda^{\mathrm T}\nabla\varepsilon_{t,
d_{1}}(\boldeta))^{2}> \epsilon/2.
\end{equation}
Denote $\bolds\lambda^{\mathrm T}(\nabla\varepsilon_{t,
d_{1}}(\boldeta)-
\nabla\varepsilon_{t}(\boldeta))$ by $R_{t}(\bolds\lambda,
\boldeta)$
and $\sigma^{-1}(\operatorname{var}(\bolds\lambda^{\mathrm T}\nabla
\varepsilon_{t, d_{1}}(\boldeta)))^{1/2}$
by $g_{t}(\sigma, \bolds\lambda,
\boldeta)$.
Since
$\bolds\lambda^{\mathrm T}\nabla\varepsilon_{t,
d_{1}}(\boldeta)/g_{t}(\sigma, \bolds\lambda,
\boldeta)$
can be written as
$\sum_{j=1}^{d_{1}} c_{j}\varepsilon_{t-j}$
with $\sum_{j=1}^{d_{1}}c_{j}^{2}=1$,
(\ref{equ3.19}) and (\ref{equ3.24}) imply that
for any $\bolds\lambda\times\boldeta\in\bolds\Lambda
\times\Pi$ and $t > d_{1}$,
%
%
\begin{eqnarray}\label{equ3.25}
&&P\bigl(s_{1} < \bolds\lambda^{\mathrm T}\nabla\varepsilon
_{t}(\boldeta) \leq s_{2}| \mathcal{F}_{t-d_{1}} \bigr)
\nonumber\\
&&\qquad= P \bigl(s_{1} + R_{t}(\bolds\lambda, \boldeta) <
\bolds\lambda^{\mathrm T}\nabla\varepsilon_{t,
d_{1}}(\boldeta) \leq s_{2} + R_{t}(\bolds\lambda, \boldeta)
| \mathcal{F}_{t-d_{1}}\bigr)
\nonumber\\
&&\qquad= P \biggl( \frac{s_{1} + R_{t}(\bolds\lambda, \boldeta)}
{g_{t}(\sigma, \bolds\lambda,
\boldeta)} <
\frac{\bolds\lambda^{\mathrm T}\nabla\varepsilon_{t,
d_{1}}(\boldeta)}
{g_{t}(\sigma, \bolds\lambda,
\boldeta)}
\leq
\frac{s_{2} + R_{t}(\bolds\lambda, \boldeta)}{
g_{t}(\sigma, \bolds\lambda,
\boldeta)
} \Big| \mathcal{F}_{t-d_{1}}\biggr) \\
&&\qquad\leq M_{1}\biggl(\frac{\sigma(s_{2} - s_{1})}{\sqrt{\epsilon
/2}}\biggr)^{\alpha_{1}}
\qquad\mbox{a.s.},\nonumber
\end{eqnarray}
provided $0 < s_{2}-s_{1} \leq(\xi\sqrt{\epsilon/2})/\sigma$.
In view of (\ref{equ3.25}), (C2) holds
with $d=d_{1}$, $M=M_{1}(\sigma\sqrt{2/\epsilon})^{\alpha_{1}}$,
$\alpha=\alpha_{1}$ and
$\delta=(\xi\sqrt{\epsilon/2})/\sigma$.

On the other hand,
it is shown in Appendix~\ref{appB} that
there exists $\tau^{**}>0$ such that
for any $\boldeta_{1}, \boldeta_{2} \in\Pi$, with $\|\boldeta
_{2}-\boldeta_{1}\|< \tau^{**}$,
%
%
\begin{equation}\label{equ3.26}
\|\nabla\varepsilon_{t}(\boldeta_{2})- \nabla\varepsilon
_{t}(\boldeta
_{1})\|
\leq\|\boldeta_{2}-\boldeta_{1}\|
\tilde{B}_{t},
\end{equation}
where $\tilde{B}_{t}$ are nonnegative random variables satisfying
%
%
\begin{equation}\label{equ3.27}
\sup_{t \geq1}{\mathrm E}(\tilde{B}^{2}_{t}) <\infty.
\end{equation}
Combining (\ref{equ3.26}) and (\ref{equ3.27}), we obtain (C3).
Finally, the proof is completed by noting that
(C4) is an immediate consequence of (\ref{equ3.26}), (\ref{equ3.27}),
(\ref{equ3.14}), (\ref{equ3.15}) and the compactness of $\Pi$.
\end{pf}
\begin{remark}\label{rem1}
In the proof of Theorem~\ref{theo3.1},
(\ref{equ3.19}) plays the same role as that of (C2)
in the proof of Theorem~\ref{theo2.1}.
When $\varepsilon_{t}$'s are normally distributed, (\ref{equ3.19})
is satisfied with $M_{1}=(2 \pi\sigma^{2})^{-1/2}$, $\alpha_{1}=1$
and any $\xi> 0$. In addition, when $\varepsilon_{t}$'s are
i.i.d. with an integrable characteristic function, (\ref{equ3.19}) is
satisfied with any $\xi> 0$, $\alpha_{1}=1$
and some $M_{1}>0$. For more details, see Lemma 4 of~\cite{11}.
An extension of Theorem~\ref{theo3.1} to autoregressive fractionally integrated
moving average models (ARFIMA)
has also been obtained by the authors.
However, since the proof of this extension is quite involved,
the details will be reported elsewhere.
\end{remark}
\begin{theorem}\label{theo3.2}
Under the same assumptions as in Theorem~\ref{theo3.1},
\textup{(C1)}--\textup{(C4)} hold for
$\mathbf{f}_{t}(\bolds\theta)=\varepsilon_{t}(\boldeta
)-\varepsilon
_{t}(\boldeta_{0})$,
$\boldTheta=\tilde{\Pi}$ and $\mathcal{F}_{t}=\sigma\{\varepsilon
_{t-1}, \varepsilon_{t-2}, \ldots\}$,
and hence by Theorem~\ref{theo2.1},
(\ref{equ3.9}) follows.
\end{theorem}

The proof of Theorem~\ref{theo3.2} is omitted, since it is similar to
the proof of Theorem~\ref{theo3.1}.
Using Theorems~\ref{theo2.2},~\ref{theo3.1} and~\ref{theo3.2}
and Lemma~\ref{lemB.1} of Appendix~\ref{appB},
the next theorem, whose proof is deferred to Appendix~\ref{appB}, establishes
moment bounds
for $n^{1/2}(\hat{\boldeta}_{n}-\boldeta_{0})$.
\begin{theorem}\label{theo3.3}
Assume that the assumptions of Theorem~\ref{theo3.1} hold
and for some $q_{1}> q \geq1$,
%
%
\begin{equation}\label{equ3.28}
\sup_{t \geq1}{\mathrm E}|\varepsilon_{t}|^{4q_{1}}< \infty.
\end{equation}
Then, (\ref{equ3.6}) follows.
\end{theorem}

As an application of Theorem~\ref{theo3.3}, an asymptotic expression for the MSPE
of~$\hat{\boldeta}_{n}$,
${\mathrm E}\{y_{n+1}-g_{n+1}(\hat{\boldeta}_{n})\}^{2}$,
is given in Theorem~\ref{theo3.4} below.
\begin{theorem}\label{theo3.4}
Assume that the assumptions of Theorem~\ref{theo3.1} hold.
Moreover, let $\varepsilon_{t}$ be i.i.d.
random variables satisfying for some $q_{1}>18$,
%
%
\begin{equation}\label{equ3.29}
{\mathrm E}|\varepsilon_{1}|^{q_{1}}< \infty.
\end{equation}
Then,
%
%
\begin{equation}\label{equ3.30}
\lim_{n \to\infty}
n[{\mathrm E}\{y_{n+1}-g_{n+1}(\hat{\boldeta}_{n})\}^{2}- \sigma
^{2}]=\bar{p}\sigma^{2}.
\end{equation}
\end{theorem}
\begin{pf}
Let $\delta_{1}$ be any positive number such that $B_{\delta
_{1}}(\bolds\eta_{0}) \subset\Pi$
and define $A_{n}=\{\hat{\bolds\eta}_{n}
\in B_{\delta_{1}}(\bolds\eta_{0})\}$
and $A^{c}_{n}=\{\hat{\bolds\eta}_{n} \in\tilde{\Pi}=
\Pi-B_{\delta_{1}}(\bolds\eta_{0})\}$.
By Taylor's theorem,
%
%
\begin{eqnarray}\label{equ3.31}
&& n^{1/2}\bigl(y_{n+1}-g_{n+1}(\hat{\bolds\eta}_{n})-\varepsilon_{n+1}\bigr)\nonumber\\
&&\hspace*{2.8pt}\qquad= n^{1/2}
(\nabla\varepsilon_{n+1}(\bolds\eta_{0}))^{\mathrm T}
(\hat{\bolds\eta}_{n}-\bolds\eta_{0})I_{A_{n}}
\nonumber\\
&&\qquad\quad{}+
\frac{n^{1/2}}{2}(\hat{\bolds\eta}_{n}-\bolds\eta_{0})^{\mathrm T}
\nabla^{2}\varepsilon_{n+1}(\bolds\eta^{*})
(\hat{\bolds\eta}_{n}-\bolds\eta_{0})I_{A_{n}}
\\
&&\qquad\quad{}+ n^{1/2}
\bigl(\varepsilon_{n+1}(\hat{\bolds\eta}_{n})-
\varepsilon_{n+1}(\bolds\eta_{0})\bigr)I_{A^{c}_{n}}\nonumber\\
&&\qquad:=
\mathrm{(I)}+\mathrm{(II)}+\mathrm{(III)},\nonumber
\end{eqnarray}
where $\|\bolds\eta^{*}-\bolds\eta_{0}\|\leq
\|\hat{\bolds\eta}_{n}-\bolds\eta_{0}\|$.
In view of (\ref{equ3.31}), (\ref{equ3.30}) holds immediately if one
can show that
%
%
\begin{eqnarray}\label{equ3.32}
\lim_{n \to\infty} {\mathrm E}\mathrm{(I)}^{2}&=&\bar{p}\sigma^{2},\\
\label{equ3.33}
\lim_{n \to\infty} {\mathrm E}(\mathrm{II})^{2}&=&0, \\
\label{equ3.34}
\lim_{n \to\infty} {\mathrm E}(\mathrm{III})^{2}&=&0.
\end{eqnarray}

By utilizing the martingale CLT (cf.~\cite{9}) and
a truncation argument in~\cite{10}, it can be shown that
%
%
\begin{equation}\label{equ3.35}
n^{1/2}\{(\nabla\varepsilon_{n+1}(\bolds\eta_{0}))^{\mathrm T}
(\hat{\bolds\eta}_{n}-\bolds\eta_{0})\}
I_{A_{n}} \Rightarrow
\mathbf{F}^{\mathrm T}\mathbf{Q},
\end{equation}
where
$\mathbf{Q}$ is distributed as $N(\mathbf{0}, \sigma^{2}\Gamma^{-1})$
with $\Gamma=\lim_{t \to\infty} {\mathrm E}\{\nabla\varepsilon
_{t}(\bolds\eta_{0})
(\nabla\varepsilon_{t}(\bolds\eta_{0}))^{\mathrm T}\}$,
and $\mathbf{F}$, satisfying ${\mathrm E}(\mathbf{F})=\mathbf{0}$
and ${\mathrm E}(\mathbf{F}\mathbf{F}^{\mathrm T})=\Gamma$,
is independent of $\mathbf{Q}$.
Let $2<r\leq18/5$. Then, it follows from H\"{o}lder's inequality,
Theorem~\ref{theo3.3}, (\ref{equ3.15}) and (\ref{equ3.29}) that
\begin{eqnarray*}
&&
\mathrm{E}\{|n^{1/2}(\nabla\varepsilon_{n+1}(\bolds\eta_{0}))^{\mathrm T}
(\hat{\bolds\eta}_{n}-\bolds\eta_{0})|^{r}\}\\
&&\qquad\leq{\mathrm E}\{
\|n^{1/2}(\hat{\bolds\eta}_{n}-\bolds\eta_{0})\|^{r}
\|\nabla\varepsilon_{n+1}(\bolds\eta_{0})\|^{r}\} \\
&&\qquad\leq\bigl({\mathrm E}
\|n^{1/2}(\hat{\bolds\eta}_{n}-\bolds\eta_{0})\|
^{5r/4}\bigr)^{4/5}
({\mathrm E}
\|\nabla\varepsilon_{n+1}(\bolds\eta_{0})\|^{5r})^{1/5}=O(1),
\end{eqnarray*}
which implies the uniform integrability of
$n\{(\nabla\varepsilon_{n+1}(\bolds\eta_{0}))^{\mathrm T}
(\hat{\bolds\eta}_{n}-\bolds\eta_{0})\}^{2}I_{A_{n}}$.
Combing this with (\ref{equ3.35}) yields
\[
\lim_{n \to\infty}
\mathrm{E}[n\{(\nabla\varepsilon_{n+1}(\bolds\eta_{0}))^{\mathrm T}
(\hat{\bolds\eta}_{n}-\bolds\eta_{0})\}^{2}I_{A_{n}}]=
\mathrm{E} (\mathbf{F}^{\mathrm T}\mathbf{Q})^{2}=\bar{p}\sigma^{2},
\]
and hence (\ref{equ3.32}) follows.
Moreover, applying Theorems~\ref{theo3.2} and~\ref{theo3.3}, (\ref{equ3.29}) and an
argument similar to that
used to prove (\ref{equB.8}) and (\ref{equB.12}) of Appendix~\ref{appB}, it is
shown in Appendix B of~\cite{4}
that (\ref{equ3.33}) and (\ref{equ3.34}) are also true.
Consequently, the desired conclusion (\ref{equ3.30}) holds.
\end{pf}
\begin{remark}\label{rem2}
Note that the moment restriction (\ref{equ3.29})
is stronger than necessary for the proofs of (\ref{equ3.32}) and (\ref
{equ3.34}).
On the other hand, since
(\ref{equ3.33}) requires that
${\mathrm E}\|n^{1/2}(\hat{\boldeta}-\boldeta_{0})\|^{q}= O(1)$ holds
with $q=9/2$
(see Appendix B of~\cite{4}),
it seems that one cannot easily weaken (\ref{equ3.29})
because Theorem~\ref{theo3.3} constitutes a key tool in verifying (\ref{equ3.33}).
\end{remark}

In the special case of $p_{2}=0$ (the pure AR case),
equation (\ref{equ3.30}) was examined by Fuller and Hasza~\cite{8},
Kunitomo and Yamamoto~\cite{13} and Ing~\cite{10}.
In addition, for the case $p_{2}>0$, equation (\ref{equ3.30})\vadjust{\goodbreak}
was also considered in Yamamoto~\cite{21}, but
a rigorous proof of (\ref{equ3.30}) is still lacking in the literature.
By establishing a set of uniform moment bounds, this paper
offers a rigorous proof of (\ref{equ3.30}) for the ARMA case.

Equation (\ref{equ3.30}) implies that
when two competing ARMA models are entertained,
the one having fewer estimated parameters also
possesses a smaller MSPE, up to terms of order $n^{-1}$.
As a result, the principle of parsimony (e.g., Tukey~\cite{20}), which
roughly asserts that mathematical models with the smallest number of
parameters are preferred,
is now endowed with a precise meaning in the context of ARMA modelling.
When $p_{2}=0$,
(\ref{equ3.30}) was established in Akaike~\cite{1} using an ad-hoc
argument, which immediately led him to
develop the final prediction error criterion,
\[
\frac{n+\bar{p}}{(n-\bar{p})n}\sum_{t=1}^{n}\bigl(y_{t}-g_{t}(\hat
{\boldeta
}_{n})\bigr)^{2}
\]
that is commonly used for AR model selection with optimal prediction
efficiency; see Shibata~\cite{18}
or Ing and Wei~\cite{12}. Under this perspective, a contribution of
(\ref{equ3.30}) is that
it provides a theoretical foundation for the construction of the FPE
criterion for ARMA models.
The issue of whether the FPE criterion (or its variants)
is asymptotically efficient (in the sense of~\cite{12} or~\cite{18})
in ARMA
model selection still remains open, however.

As a final remark, we note that
(\ref{equ3.30}) is obtained based on Theorems~\ref{theo2.1} and~\ref{theo2.2}.
Moreover, since these theorems provide a useful device for exploring
the moment properties of least squares estimates in (nonlinear)
stochastic regression models, their applications to prediction or model
selection
in models beyond the ARMA case are anticipated.

\begin{appendix}\label{app}
\section{\texorpdfstring{Proofs of (\lowercase{\protect\ref{equ2.3}}) and (\lowercase{\protect\ref{equ2.26}})}
{Proofs of (2.3) and (2.26)}}\label{appA}

\vspace*{-8pt}

\begin{pf*}{Proof of (\ref{equ2.3})}
Let $m = \lfloor\{l_{1}(r+2k)+r+k+2q\}/\alpha\rfloor+1$ with $l_{1}>q$.
We only prove (\ref{equ2.3}) for the case of $l=0$ and $j=1$
since the other\vspace*{1pt} cases can be similarly verified.
First, define
$A(u)= \{\sum_{i=0}^{m-1}\sup_{\bolds\theta\in
\boldTheta}
\|\mathbf{f}_{(i+1)d+1}(\bolds\theta)\|^{2} \leq u^{l_{1}/q}/r\}$ and
$B(u)= \{\sum_{i=0}^{m-1} B_{(i+1)d+1} \leq u^{l_{1}/q}/k^{1/2}\}$,
where $B_{t}$ are random variables defined in (C3).
Then, the left-hand side of (\ref{equ2.3}) (with $l=0$ and $j=1$)
is bounded by
%
%
\begin{eqnarray}\label{equA.1}\hspace*{-3pt}
&& K_{0}+\int_{K_{0}}^{\infty} P\Biggl\{\sup_{\bolds\theta\in
\boldTheta}
\Biggl(\inf_{\|\mathbf{y}\|=1}\sum_{i=0}^{m-1}
\bigl(\mathbf{y}^{\mathrm T}\mathbf{f}_{(i+1)d+1}(\bolds\theta)\bigr)^{2}\Biggr)^
{-q} > u\Biggr\} \,du \nonumber\\
&&\qquad= K_{0} +\int_{K_{0}}^{\infty} P\Biggl\{
\inf_{\bolds\theta\in\boldTheta}
\inf_{\|\mathbf{y}\|=1}\sum_{i=0}^{m-1}
\bigl(\mathbf{y}^{\mathrm T}\mathbf{f}_{(i+1)d+1}(\bolds\theta)\bigr)^{2}
< u^{-1/q}\Biggr\} \,du \nonumber\\
&&\qquad\leq
K_{0}+ \int_{K_{0}}^{\infty}\! P\Biggl\{\inf_{\bolds\theta\in
\boldTheta}
\inf_{\|\mathbf{y}\|=1}\!\sum_{i=0}^{m-1}
\bigl(\mathbf{y}^{\mathrm T}\mathbf{f}_{(i+1)d+1}(\bolds\theta)\bigr)^{2}
< u^{-1/q}, A(u), B(u)\!\Biggr\} \,du \hspace*{-8pt}\\
&&\qquad\quad{} + \int_{K_{0}}^{\infty} P(A^{c}(u)) \,du +
\int_{K_{0}}^{\infty} P(B^{c}(u)) \,du \nonumber\\
&&\qquad\equiv K_{0} + \mathrm{(I)} + \mathrm{(II)} + \mathrm{(III)},
\nonumber
\end{eqnarray}
where $K_{0}=K_{0}(l_{1}, \delta, q, k,\tau)$ is a positive number to be
specified later
and $A^{c}(u)$ and $B^{c}(u)$
denote the complements of $A(u)$ and $B(u)$, respectively.
Since $l_{1}>q$, by (C3), (C4)
and Chebyshev's inequality,
it follows that for $n$ large,
%
%
\begin{equation}\label{equA.2}
\mathrm{(II)} \leq C_{1}^{\ast} \quad\mbox{and}\quad \mathrm{(III)} \leq C_{2}^{\ast},
\end{equation}
where $C_{1}^{\ast}$ and $C_{2}^{\ast}$ are some positive constants
depending on $C_{1}, C_{2}, \alpha, l_{1}, r, k, q$
and $K_{0}$.

To deal\vspace*{1pt} with (I), consider the hypersphere $\mathbf{S}_{r} = \{
\mathbf{y}\dvtx \mathbf{y} \in R^{r}, \| \mathbf{y}\| = 1 \}$ and
the hypercube $\mathbf{H}^{r}(u) = [ 1 - 2 u^{- ( l_{1} + 1)/2q}
( \lfloor u^{( l_{1} + 1)/2q} \rfloor+ 1 ), 1]^{r}, u > 0$. Note
first that $\mathbf{S}_{r} \subseteq\mathbf{H}^{r}(u)$ for any
$u>0$. Divide $\mathbf{H}^{r}(u)$ into sub-hypercubes of equal size,
each of which has an edge length of $2 u^{- ( l_{1} + 1)/2q}$ and a
circumscribed circle of radius $\sqrt{r} u^{- ( l_{1} + 1)/2q}$.
Denote these\vspace*{1pt} sub-hypercubes by $\tilde{B}_{i}(u)$, $1 \leq i \leq
m^{\ast}=( \lfloor u^{( l_{1} + 1)/2q} \rfloor+ 1 )^{r}$. Letting
$G_{i}(u)=\mathbf{S}_{r} \cap\tilde{B}_{i}(u)$ and $\{G_{v_{i}}(u),
i= 1, \ldots, m^{\ast\ast}\}$ denote the collection of nonempty
$G_{i}(u)$'s, it follows that $\mathbf{S}_{r} = \bigcup_{i=1}^{m^{\ast
\ast}} G_{v_{i}}(u)$ with $m^{\ast\ast} \leq( \lfloor u^{( l_{1} +
1)/2q} \rfloor+ 1 )^{r}$. On the other hand, since $\boldTheta$ is
a bounded subset in $R^{k}$, there is a positive integer $g$ such
that for any $u>0$, $\boldTheta\subseteq\mathbf{H}_{g}^{k} (u)=
[g-2gu^{-(l_{1}+1/2)/q}(\lfloor u^{(l_{1}+1/2)/q}\rfloor+1),
g]^{k}$. We can\vspace*{-2pt} similarly divide $\mathbf{H}_{g}^{k}(u)$ into
equal-sized sub-hypercubes $\tilde{W}_{i}(u), i=1, \ldots, e^{*}$,
where the edge length of $\tilde{W}_{i}(u)$ is $2 u^{- ( l_{1} + 1/2
) q^{-1}}$ and $e^{\ast}= g^{k} ( \lfloor u^{( l_{1} + 1/2 )/ q}
\rfloor+ 1 )^{k}$. In addition, it holds that $\boldTheta=
\bigcup_{i=1}^{e^{\ast\ast}}J_{v_{i}}(u)$, where with
$J_{i}(u)=\boldTheta\cap\tilde{W}_{i}(u)$, $\{J_{v_{i}}(u), i= 1,
\ldots, e^{\ast\ast}\}$ denotes the collection of nonempty
$J_{i}(u)$'s. By observing
\begin{eqnarray*}
&& \Biggl\{\inf_{\bolds\theta\in\boldTheta}
\inf_{\|\mathbf{y}\|=1} \sum_{i=0}^{m-1} \bigl(\mathbf{y}^{\mathrm T}
\mathbf{f}_{(i+1)d+1}(\bolds\theta)\bigr)^{2}
< u^{-1/q}\Biggr\} \\
&&\qquad=
\bigcup_{s=1}^{e^{\ast\ast}}\bigcup_{j=1}^{m^{\ast\ast}}
\Biggl\{\inf_{\bolds\theta\in J_{v_{s}}(u)} \inf_{\mathbf{y}
\in G_{v_{j}}(u)} \sum_{i=0}^{m-1} \bigl(\mathbf{y}^{\mathrm T}
\mathbf{f}_{(i+1)d+1}(\bolds\theta)\bigr)^{2}
< u^{-1/q}\Biggr\},
\end{eqnarray*}
one has
%
%
\begin{eqnarray}\label{equA.3}
&& P\Biggl(\inf_{\bolds\theta\in\boldTheta}
\inf_{\|\mathbf{y}\|=1} \sum_{i=0}^{m-1} \bigl(\mathbf{y}^{\mathrm
T}\mathbf
{f}_{(i+1)d+1}(\bolds\theta)\bigr)^{2}
< u^{-1/q}, A(u), B(u)\Biggr) \nonumber\\[-8pt]\\[-8pt]
&&\qquad \leq
\sum_{s=1}^{e^{\ast\ast}}\sum_{j=1}^{m^{\ast\ast}} P\Biggl(
\bigcap_{i=0}^{m-1} C^{(s,j)}_{i}(u)
\Biggr),
\nonumber
\end{eqnarray}
where
\begin{eqnarray*}
C^{(s,j)}_{i}(u)
&=& \biggl\{\inf_{\bolds\theta\in J_{v_{s}}(u)}
\inf_{\mathbf{y} \in G_{v_{j}}(u)}
\bigl|\mathbf{y}^{\mathrm T}\mathbf{f}_{(i+1)d+1}(\bolds\theta)\bigr| <
u^{{-1}/({2q})},\\
&&\hspace*{5.7pt} B_{(i+1)d+1} \leq\frac{u^{l_{1}/q}}{k^{1/2}},
\sup_{\bolds\theta\in\boldTheta}
\bigl\|\mathbf{f}_{(i+1)d+1}(\bolds\theta)\bigr\| \leq\frac{u^{l_{1}/2q}}
{r^{1/2}} \biggr\}.
\end{eqnarray*}
Let $\mathbf{y}_{j} \in G_{v_{j}}(u), j=1, \ldots, m^{\ast\ast}$, and
$\bolds\theta_{s} \in J_{v_{s}}(u),
s=1, \ldots, e^{\ast\ast}$,
be arbitrarily chosen.
Then, for any $\mathbf{y} \in G_{v_{j}}(u)$ and
$\bolds\theta\in J_{v_{s}}(u)$,
\begin{eqnarray*}
\bigl|\mathbf{y}^{\mathrm T}_{j}\mathbf{f}_{(i+1)d+1}(\bolds\theta_{s})\bigr|
&\leq&\|\mathbf{y}_{j}-\mathbf{y}\| \bigl\|\mathbf{f}_{(i+1)d+1}
(\bolds\theta_{s})\bigr\| \\
&&{}+ \|\mathbf{y}\|\bigl\|\mathbf{f}_{(i+1)d+1}(\bolds\theta_{s})-
\mathbf{f}_{(i+1)d+1}(\bolds\theta)\bigr\|\\
&&{} + \bigl|\mathbf{y}^{\mathrm T}\mathbf{f}_{(i+1)d+1}(\bolds\theta)\bigr|.
\end{eqnarray*}
Combining this with (C3) yields that on the set $C^{(s, j)}_{i}(u)$
with $u > (2 k^{1/2}/\break\tau)^{q/(l_{1}+1/2)}$,
\begin{eqnarray*}
&& \bigl|\mathbf{y}^{\mathrm
T}_{j}\mathbf{f}_{(i+1)d+1}(\bolds\theta_{s})\bigr|\\
&&\qquad\leq2 \sqrt{r} u^{-(l_{1}+1)/2q} \sup_{\bolds\theta\in
\boldTheta}
\bigl\|\mathbf{f}_{(i+1)d+1}(\bolds\theta)\bigr\| \\
&&\qquad\quad{} + 2 \sqrt{k} u^{-(l_{1}+1/2)/q} B_{(i+1)d+1} +
\inf_{\bolds\theta\in J_{v_{s}}(u)} \inf_{\mathbf{y} \in
G_{v_{j}}(u)} \bigl|\mathbf{y}^{\mathrm T}\mathbf{f}_{(i+1)d+1}(\bolds
\theta)\bigr|\\
&&\qquad\leq5 u^{-1/2q}
\end{eqnarray*}
and hence
%
%
\begin{equation}\label{equA.4}
C^{(s,j)}_{i}(u) \subseteq D^{(s,j)}_{i} (u) := \bigl\{\bigl|\mathbf{y}^{\mathrm T}_{j}
\mathbf{f}_{(i+1)d+1}(\bolds\theta_{s})\bigr|
\leq5 u^{-1/2q}\bigr\}.
\end{equation}
In view of (\ref{equA.3}) and (\ref{equA.4}), it follows that for $u >
(2 k^{1/2}/\tau)^{q/(l_{1}+1/2)}$,
%
%
\begin{eqnarray}\label{equA.5}
&& P\Biggl(\inf_{\bolds\theta\in\boldTheta}
\inf_{\|\mathbf{y}\|=1} \sum_{i=0}^{m-1} \bigl(\mathbf{y}^{\mathrm
T}\mathbf
{f}_{(i+1)d+1}(\bolds\theta)\bigr)^{2}
< u^{-1/q}, A(u), B(u)\Biggr) \nonumber\\[-8pt]\\[-8pt]
&&\qquad\leq
\sum_{s=1}^{e^{\ast\ast}}\sum_{j=1}^{m^{\ast\ast}} P\Biggl(
\bigcap_{i=0}^{m-1} D^{(s,j)}_{i}(u)
\Biggr).
\nonumber
\end{eqnarray}
Observe that
\[
P\Biggl(
\bigcap_{i=0}^{m-1} D^{(s,j)}_{i}(u)
\Biggr) = \mathrm{E}\Biggl\{ \prod_{i=0}^{m-2} I_{D^{(s,j)}_{i}(u)}
P\bigl(D^{(s,j)}_{m-1}(u)|\mathcal{F}_{(m-1)d+1}\bigr)
\Biggr\},
\]
where $I_{D^{(s,j)}_{i}(u)}$ denotes the indicator function of the set
$D^{(s,j)}_{i}(u)$.
This, together with (C2), implies that
for $u > (10/\delta)^{2q}$,
all $1 \leq s \leq e^{\ast\ast}$, all $1 \leq j \leq m^{\ast\ast}$
and $n$ large,
\[
P\Biggl(
\bigcap_{i=0}^{m-1} D^{(s,j)}_{i}(u)
\Biggr) \leq M (10)^{\alpha} u^{-\alpha/2q} \mathrm{E}\Biggl\{
\prod_{i=0}^{m-2} I_{D^{(s,j)}_{i}(u)} \Biggr\}.
\]
Repeating the same argument $m-1$ times, one has
%
%
\begin{equation}\label{equA.6}
P\Biggl(
\bigcap_{i=0}^{m-1} D^{(s,j)}_{i}(u)
\Biggr)
\leq M^{m}(10)^{ m \alpha} u^{- m \alpha/2q}.
\end{equation}
Taking $K_{0} > \max\{(10/\delta)^{2q}, (2k^{1/2}/\tau
)^{q/(l_{1}+1/2)}, 1\}$,
it follows from (\ref{equA.5}), (\ref{equA.6}) and $m > \{
l_{1}(r+2k)+r+k+2q\}/\alpha$ that
%
%
\begin{eqnarray}\label{equA.7}
\mathrm{(I)} &\leq&\int_{K_{0}}^{\infty} \sum_{s=1}^{e^{\ast\ast}}
\sum_{j=1}^{m^{\ast\ast}} P\Biggl(
\bigcap_{i=0}^{m-1} D^{(s,j)}_{i}(u) \Biggr) \,du \nonumber\\
&\leq& 2^{r+k} g^{k} M^{m} (10)^{\alpha m}
\int_{K_{0}}^{\infty} u^{-\{1/({2q})\}\{\alpha m -(l_{1}+1)r-
(2l_{1}+1)k\}} \,du \\
&=& 2^{r+k} g^{k} M^{m} (10)^{\alpha m} \{C(q, \alpha, m, l_{1}, r, k)\}^{-1}
K_{0}^{-C(q, \alpha, m, l_{1}, r, k)},
\nonumber
\end{eqnarray}
where $C(q, \alpha, m, l_{1}, r, k) =
\{\alpha m - (l_{1}+1)r - (2l_{1}+1)k - 2q\}/2q$.
Consequently, (\ref{equ2.3}) is ensured by (\ref{equA.1}), (\ref
{equA.2}) and (\ref{equA.7}).
\end{pf*}
\begin{pf*}{Proof of (\ref{equ2.26})}
Let $C^{*}_{4}=\delta^{*^{q}}_{1}k^{q}\bar{M}^{2q}$.
Since $\delta^{*}_{1}$, defined at the beginning
of the proof of Theorem~\ref{theo2.2}, is smaller than $3^{-1}k^{-1}\bar{M}^{-2}$,
it follows that $C^{*}_{4}<3^{-q}$.
By the Cauchy--Schwarz inequality and (\ref{equ2.12})--(\ref{equ2.14}),
one has
%
%
\begin{eqnarray}\label{equA.8}
{\mathrm E}\mathrm{(III)} &\leq& \delta_{1}^{*^{2q}}k^{q}n^{q/2}
{\mathrm E}\bigl(R^{q}_{n}R^{q/2}_{1n}R^{q/2}_{2n}
I_{\{R_{n}R^{1/2}_{1n}R^{1/2}_{2n}>\bar{M}^{2}\}}\bigr) \nonumber\\
&&{} + C^{*}_{4} {\mathrm E}\bigl(\|n^{1/2}(\hat{\bolds\theta}_{n}-
\bolds\theta_{0})\|^{q}I_{A_{n}}\bigr)\nonumber\\
&\leq&\delta_{1}^{*^{2q}}k^{q}n^{q/2}
\{{\mathrm E}(R^{2q}_{n}R^{q}_{1n}R^{q}_{2n})\}^{1/2}
\nonumber\\[-8pt]\\[-8pt]
&&{}\times\{P(R_{n}>
\bar{M})+P(R_{1n}> \bar{M})+P(R_{2n}> \bar{M})\}^{1/2} \nonumber\\
&&{}+
C^{*}_{4} {\mathrm E}\bigl(\|n^{1/2}(\hat{\bolds\theta}_{n}
-\bolds\theta_{0})\|^{q}I_{A_{n}}\bigr) \nonumber\\
&=& O(1)\{{\mathrm E}(R^{2q}_{n}R^{q}_{1n}R^{q}_{2n})\}^{1/2}+C^{*}_{4}
{\mathrm E}\bigl(\|n^{1/2}(\hat{\bolds\theta}_{n}
-\bolds\theta_{0})\|^{q}I_{A_{n}}\bigr).
\nonumber
\end{eqnarray}
In addition, ${\mathrm E}(R^{2q}_{n}R^{q}_{1n}R^{q}_{2n})=O(1)$
follows from H\"{o}lder's inequality, (\ref{equ2.9}), (\ref{equ2.10})
and (\ref{equ2.19}).
Combining this with (\ref{equA.8}) yields (\ref{equ2.26}).
\end{pf*}

\section{\texorpdfstring{Proofs of (\lowercase{\protect\ref{equ3.26}}), (\lowercase{\protect\ref{equ3.27}}) and Theorem \lowercase{\protect\ref{theo3.3}}}
{Proofs of (3.26), (3.27) and Theorem 3.3}}\label{appB}

Throughout this Appendix,
$\mathbf{J}(m,\bar{p}), 1 \leq m \leq\bar{p}$,
denotes the set
$\{(j_{1}, \ldots,\break j_{m})\dvtx j_{1}<\cdots<j_{m}, j_{i}\in\{1, \ldots,
\bar{p}\}
\mbox{ for } 1\leq i \leq m\}$,
and for
$\mathbf{j}=(j_{1}, \ldots, j_{m}) \in\mathbf{J}(m,\bar{p})$ and\vadjust{\goodbreak}
smooth function $w=w(\boldxi)=w(\xi_{1}, \ldots, \xi_{\bar{p}})$,
$\mathbf{D}_{\mathbf{j}} w$ denotes
the partial derivative
$\partial^{m} w/\partial\xi_{j_{1}}, \ldots, \partial\xi_{j_{m}}$.
Before proving (\ref{equ3.26}) and (\ref{equ3.27}), we note that
according to
(\ref{equ3.3})--(\ref{equ3.5}), (\ref{equ3.10})--(\ref{equ3.14}) and
the compactness of $\Pi$,
$
(\nabla^{2} \varepsilon_{t}(\boldeta))_{i,j}
=\sum_{s=1}^{t-2}c_{s, ij}(\boldeta)\varepsilon_{t-1-s}$,
where $c_{s, ij}(\boldeta)$ are continuously differentiable on $\Pi$
and satisfy, for some $D_{1}, D_{2}>0$ (independent of $i, j$ and $s$),
%
%
\begin{equation}\label{equB.1}
\sup_{\boldeta\in\Pi}
|c_{s, ij}(\boldeta)| \leq D_{1}\operatorname{exp}(-D_{2}s).
\end{equation}
Moreover, there exists a small positive number $\tau^{*}$ such that
%
%
\begin{eqnarray}\label{equB.2}
\sup_{\boldeta\in\Pi^{*}}
|\mathbf{D}_{\mathbf{j}} b_{s}(\boldeta)| &\leq& D_{3}
\operatorname{exp}(-D_{6}s),\\
\label{equB.3}
\max_{\mathbf{j} \in\mathbf{J}(m, \bar{p}), 1\leq m \leq\bar
{p}}\sup
_{\boldeta\in\Pi^{*}}
\bigl|\mathbf{D}_{\mathbf{j}} b^{(l)}_{s}(\boldeta)\bigr| &\leq& D_{4}
\operatorname{exp}(-D_{6}s), \\
\label{equB.4}
\max_{\mathbf{j} \in\mathbf{J}(m, \bar{p}), 1\leq m \leq\bar
{p}}\sup
_{\boldeta\in\Pi^{*}}
|\mathbf{D}_{\mathbf{j}} c_{s, ij}(\boldeta)| &\leq& D_{5}
\operatorname{exp}(-D_{6}s),
\end{eqnarray}
where $\Pi^{*}=\bigcup_{\boldeta\in\Pi}B_{\tau^{*}}(\boldeta)$ and
$D_{3}, \ldots, D_{6}$ are some positive constants independent of $i,
j, l$ and $s$.
\begin{pf*}{Proofs of (\ref{equ3.26}) and (\ref{equ3.27})}
Let $\tau^{**}=\tau^{*}/2$.
For $\|\boldeta_{2}-\boldeta_{1}\|<\tau^{**}$, it follows from
the mean value theorem for vector-valued functions that
$
\|\nabla\varepsilon_{t}(\boldeta_{2})-\nabla\varepsilon
_{t}(\boldeta
_{1})\|^{2}
\leq\|\boldeta_{2}-\boldeta_{1}\|^{2}
\|\int_{0}^{1}\nabla^{2}\varepsilon_{t}(\boldeta_{1}+v(\boldeta
_{2}-\boldeta_{1}))\,dv\|^{2}
\leq\|\boldeta_{2}-\boldeta_{1}\|^{2}(\tilde{B}_{t})^{2}$,
where
$\tilde{B}_{t}= \{
\sum_{1\leq i,j\leq\bar{p}}
\sup_{\boldeta\in\Pi^{**}}$
$(\nabla^{2}\varepsilon_{t}(\boldeta))^{2}_{i,j}
\}^{1/2}$,
with $\Pi^{**}=\break\bigcup_{\boldeta\in\Pi} B_{\tau^{**}}(\boldeta)$.
Denoting by $\bar{\Pi}^{**}$ the compact closure of $\Pi^{**}$,
one has $\bar{\Pi}^{**} \subset\Pi^{*}$, which further yields $\bar
{\Pi
}^{**} \subset\bigcup_{r=1}^{\bar{r}} B_{\tau^{*}}(\bolds\theta_{r})$,
for some $1\leq\bar{r}<\infty$ and $\bolds\theta_{1}, \ldots,\break
\bolds\theta_{\bar{r}}$ $\in\Pi$.
Hence,
$\mathrm E (\tilde{B}^{2}_{t}) \leq
\sum_{1\leq i,j\leq\bar{p}} \sum_{r=1}^{\bar{r}}
\mathrm E\{\sup_{\boldeta\in B_{\tau^{*}}(\bolds\theta
_{r})}(\nabla
^{2}\varepsilon_{t}(\boldeta))^{2}_{i,j}
\}$. Moreover,\vspace*{1pt} it follows from (B.1), (B.4)
and (3.10) of Lai~\cite{14} that
for all $1\leq i, j \leq\bar{p}$, $1\leq r \leq\bar{r}$ and $t \geq3$,
$
\mathrm E\{\sup_{\boldeta\in B_{\tau^{*}}(\bolds\theta
_{r})}(\nabla
^{2}\varepsilon_{t}(\boldeta))^{2}_{i,j}
\}< C
\sum_{s=1}^{\infty}\{\operatorname{exp}(-2D_{2}s)+\operatorname{exp}(-2D_{6}s)\}$\vspace*{1pt}
for some $C>0$ (see Appendix B of~\cite{4} for more details).
Consequently, (\ref{equ3.26}) and (\ref{equ3.27}) follow.
\end{pf*}

The next lemma, Lemma~\ref{lemB.1},
provides moment bounds for the supremums
of some random functions associated with
(\ref{equ2.8}) and (\ref{equ2.11})--(\ref{equ2.14}).
Lemma~\ref{lemB.1}, together with Theorems~\ref{theo3.1} and~\ref{theo3.2}, constitutes the major
tools for proving Theorem~\ref{theo3.3}.
\begin{lemma}\label{lemB.1}
Let $\bolds\theta_{a}$ be some point in $R^{k}, k \geq1$,
and $\delta_{1}$ be some positive number.
For $t \geq2$, define
$K_{t}(\bolds\theta)=\sum_{i=1}^{t-1}c_{i}(\bolds\theta)\epsilon_{t-i}$
and $Q_{t}(\bolds\theta)=\sum_{i=1}^{t-1}d_{i}(\bolds\theta
)\epsilon_{t-i}$,
where $\epsilon_{i}$ are independent random variables with $\mathrm{E}(\epsilon_{i})=0$
and $\mathrm{E}(\epsilon^{2}_{i})=\sigma^{2}_{\epsilon}>0$ for all $i\geq1$,
and
$c_{i}(\bolds\theta)$ and $d_{i}(\bolds\theta)$ are real-valued functions
on $B_{\delta_{1}}(\bolds\theta_{a})$.
Assume that for any $i \geq1$,
$\mathbf{j} \in\mathbf{J}(m, k)$ and $1\leq m \leq k$,
$\mathbf{D}_{\mathbf{j}}c_{i}(\bolds\theta)$
are continuous on $B_{\delta_{1}}(\bolds\theta_{a})$, and
for some $q_{1} \geq2$,
$\sup_{i \geq1}{\mathrm E}|\epsilon_{i}|^{q_{1}}< \infty$.
Then, there exists $C>0$ such that for all $n \geq2$,
%
%
\begin{eqnarray}\label{equB.5}
&&{\mathrm E}\Biggl(\sup_{\bolds\theta\in
B_{\delta_{1}}(\bolds\theta_{a})}
\Biggl|\sum_{t=2}^{n}K_{t}(\bolds\theta)\epsilon_{t}\Biggr|^{q_{1}}\Biggr) \nonumber\\
&&\qquad\leq C
n^{{q_{1}}/{2}}\Biggl[
\Biggl\{\sum_{i=1}^{n-1}
c^{2}_{i}(\bolds\theta_{a})
\Biggr\}^{{q_{1}}/{2}} \\
&&\qquad\quad\hspace*{36.5pt}{} +
\Biggl\{\sum_{i=1}^{n-1}
\max_{\mathbf{j} \in\mathbf{J}(m, k), 1\leq m \leq k}
\sup_{\bolds\theta\in
B_{\delta_{1}}(\bolds\theta_{a})}(
\mathbf{D}_{\mathbf{j}}c_{i}(\bolds\theta)
)^{2}\Biggr\}^{{q_{1}}/{2}}\Biggr].
\nonumber
\end{eqnarray}
Moreover, if
for any $i, j \geq1$,
$\mathbf{j} \in\mathbf{J}(m, k)$ and $1\leq m \leq k$,
$\mathbf{D}_{\mathbf{j}}\{c_{i}(\bolds\theta)d_{j}(\bolds\theta)\}$
are continuous on $B_{\delta_{1}}(\bolds\theta_{0})$, and for some $q_{1}
\geq2$,
$\sup_{i \geq1}{\mathrm E}|\epsilon_{i}|^{2q_{1}}< \infty$,
then there exists $C>0$ such that for all $n \geq3$,
%
%
\begin{eqnarray}\label{equB.6}
&& {\mathrm E}\Biggl(\sup_{\bolds\theta\in
B_{\delta_{1}}(\bolds\theta_{a})}
\Biggl|\sum_{t=2}^{n}K_{t}(\bolds\theta)Q_{t}(\bolds\theta
)-\mathrm{E}(K_{t}(\bolds\theta
)Q_{t}(\bolds\theta))\Biggr|^{q_{1}}\Biggr)\nonumber\\
&&\qquad \leq C \Biggl[
\Biggl\{\sum_{j=1}^{n-1}\Biggl(\sum_{l=1}^{n-j}
S_{l,l}
\Biggr)^{2} + \sum_{j=1}^{n-1} \Biggl(\sum_{l=1}^{n-j} V_{l,l}\Biggr)^{2}
\Biggr\}^{{q_{1}}/{2}}\nonumber\\
&&\qquad\quad\hspace*{14pt}{}  +
n^{({q_{1}-2})/{2}}
\sum_{j=2}^{n-1} \Biggl\{ \Biggl(
\sum_{i=1}^{j-1} \Biggl(\sum_{l=1}^{n-j}
S_{l+j-i, l}
\Biggr)^{2} \Biggr)^{{q_{1}}/{2}}\nonumber\\[-8pt]\\[-8pt]
&&\qquad\quad\hspace*{14pt}\hspace*{75.6pt}{} + \Biggl(
\sum_{i=1}^{j-1} \Biggl(\sum_{l=1}^{n-j} S_{l, l+j-i}
\Biggr)^{2} \Biggr)^{{q_{1}}/{2}}\nonumber\\
&&\qquad\quad\hspace*{14pt}\hspace*{75.6pt}{}  +
\Biggl(\sum_{i=1}^{j-1}\Biggl(\sum_{l=1}^{n-j}
V_{l+j-i, l}\Biggr)^{2} \Biggr)^{{q_{1}}/{2}}\nonumber\\
&&\qquad\quad\hspace*{14pt}\hspace*{75.6pt}{} +
\Biggl(\sum_{i=1}^{j-1}\Biggl(\sum_{l=1}^{n-j}
V_{l, l+j-i}\Biggr)^{2} \Biggr)^{{q_{1}}/{2}}
\Biggr\}
\Biggr],
\nonumber
\end{eqnarray}
where
\[
V_{i,j}=|c_{i}(\bolds\theta_{a})d_{j}(\bolds\theta_{a})|
\quad\mbox{and}\quad
S_{i, j}=\max_{\mathbf{j} \in\mathbf{J}(m, k), 1\leq m \leq k}\sup
_{\bolds\theta\in
B_{\delta_{1}}(\bolds\theta_{a})}
|\mathbf{D}_{\mathbf{j}}\{c_{i}(\bolds\theta)d_{j}(\bolds\theta)\}|.
\]
\end{lemma}

The proof of (\ref{equB.5}), given in Appendix B of~\cite{4},
is based on (3.8) of~\cite{14} and Lemma 2 of~\cite{22}.
Assuming that $\sup_{i \geq1}{\mathrm E}|\varepsilon_{i}|^{q_{1}}<
\infty$
for some $q_{1}>\max\{q, 2\}$ with $q \geq1$, (\ref{equB.5}) can be
used to justify (\ref{equ2.8}) for the ARMA case.
More precisely, applying
(\ref{equB.5}) with $K_{t}(\bolds\theta)=(\nabla^{2} \varepsilon
_{t}(\boldeta))_{i,j}$
and $\epsilon_{t}=\varepsilon_{t}$,
in conjunction with
(\ref{equB.1}) and (\ref{equB.4}), it follows that for any $\delta
_{1}>0$ with $B_{\delta_{1}}(\boldeta_{0}) \subset\Pi$,
%
%
\begin{equation}\label{equB.7}
\max_{1\leq i, j \leq\bar{p}}\mathrm{E}\Biggl(
\sup_{\boldeta\in
B_{\delta_{1}}(\boldeta_{0})}
\Biggl|n^{-1/2}\sum_{t=1}^{n}\varepsilon_{t}(\nabla^{2} \varepsilon_{t}
(\boldeta))_{i,j}\Biggr|^{q_{1}}\Biggr)= O(1).
\end{equation}
In addition, by making use of (\ref{equB.5}) with
$K_{t}(\bolds\theta)=\varepsilon_{t}(\boldeta)-\varepsilon
_{t}(\boldeta_{0})$
and $\epsilon_{t}=\varepsilon_{t}$,
the compactness of $\tilde{\Pi}$, (\ref{equ3.17}) and (\ref{equB.2}),
we obtain
%
%
\begin{equation}\label{equB.8}
\mathrm{E}\Biggl(\sup_{\boldeta\in\tilde{\Pi}}
\Biggl|n^{-1}\sum_{t=1}^{n}\varepsilon_{t}\bigl(\varepsilon_{t}(\boldeta)-
\varepsilon_{t}(\boldeta_{0})\bigr)\Biggr|^{q_{1}}\Biggr) =O(n^{- q_{1}/2}),
\end{equation}
which gives (\ref{equ2.11}) (with $q_{2}=q_{1}$ and $\nu=1/2$) for the
ARMA case.

On the other hand, (\ref{equB.6}), whose proof is also given in
Appendix B of~\cite{4},
can be viewed as a uniform version of
the first moment bound theorem of~\cite{6} and plays a key role in
verifying (\ref{equ2.12})--(\ref{equ2.14}) for the ARMA case.
Let $\bar{M}_{3}$ be any positive number larger than
$2D^{2}_{1}\sigma^{2}\sum_{l=1}^{\infty}\operatorname{exp}(-2D_{2}l)$ and
$\delta_{1}$ be any positive number
satisfying $B_{\delta_{1}}(\boldeta_{0}) \subset\Pi$,
noting that $D_{1}$ and $D_{2}$ are defined in (\ref{equB.1}).
Assume $\sup_{i \geq1}{\mathrm E}|\varepsilon_{i}|^{2q_{1}}< \infty$
for some $q_{1} \geq2q$ with $q \geq1$.
Then, by (\ref{equB.6}) with $K_{t}(\bolds\theta)=Q_{t}(\bolds
\theta
)=(\nabla^{2} \varepsilon_{t}(\boldeta))_{i,j}$
and $\epsilon_{t}=\varepsilon_{t}$, (\ref{equB.1}), (\ref{equB.4}) and
Chebyshev's inequality,
one has for any $1\leq i, j \leq\bar{p}$,
%
%
\begin{eqnarray}\label{equB.9}
&& P\Biggl(\sup_{\boldeta\in
B_{\delta_{1}}(\boldeta_{0})} n^{-1}
\sum_{t=1}^{n}(\nabla^{2} \varepsilon_{t}(\boldeta))_{i,j}^{2}
>\bar{M}_{3}
\Biggr) \nonumber\hspace*{-15pt}\\
&&\qquad\leq
P\Biggl(\sup_{\boldeta\in
B_{\delta_{1}}(\boldeta_{0})} \Biggl| n^{-1}
\sum_{t=1}^{n}
[(\nabla^{2} \varepsilon_{t}(\boldeta))_{i,j}^{2}-
\mathrm{E}\{(\nabla^{2} \varepsilon_{t}(\boldeta))_{i,j}^{2}\}
]\Biggr|^{q_{1}}
>(\bar{M}_{3}/2)^{q_{1}}
\Biggr)\hspace*{-15pt} \\
&&\qquad= O(n^{-q_{1}/2})=O(n^{-q}),
\nonumber\hspace*{-15pt}
\end{eqnarray}
which is (\ref{equ2.14}) for the ARMA case.
In addition, (\ref{equ2.12}) and (\ref{equ2.13}) for the ARMA case,
that is, for some $\bar{M}_{1}, \bar{M}_{2}>0$,
%
%
\begin{eqnarray}\quad
\label{equB.10}
P\Biggl(\sup_{\boldeta\in B_{\delta_{1}}
(\boldeta_{0})}
\lambda^{-1}_{\min}\Biggl(n^{-1}
\sum_{t=1}^{n}\nabla\varepsilon_{t}(\boldeta) (\nabla\varepsilon_{t}
(\boldeta))^{\mathrm T}\Biggr)
>\bar{M}_{1}
\Biggr) &=& O(n^{-q}), \\
\label{equB.11}
P\Biggl(\sup_{\boldeta\in B_{\delta_{1}}
(\boldeta_{0})} n^{-1}
\sum_{t=1}^{n}\|\nabla\varepsilon_{t}(\boldeta)\|^{2}
>\bar{M}_{2}
\Biggr) &=& O(n^{-q}),
\end{eqnarray}
can also be similarly verified.
With the help of these results, we are now in a position to prove
Theorem~\ref{theo3.3}.
\begin{pf*}{Proof of Theorem~\ref{theo3.3}}
Since (\ref{equ3.28}) is assumed, (\ref{equB.7})--(\ref{equB.11}) follow.
In view of Theorems~\ref{theo2.2},~\ref{theo3.1} and~\ref{theo3.2},
it remains to show that for some $q_{1}> q \geq1$ and
some small positive number $\delta_{1}$
with $B_{\delta_{1}}(\boldeta) \subset\Pi$,
%
%
\begin{equation}\label{equB.12}
\max_{1\leq i, j \leq\bar{p}, 1 \leq t \leq n}\mathrm{E}\Bigl(\sup_
{\boldeta\in B_{\delta_{1}}(\boldeta_{0})}
|(\nabla^{2} \varepsilon_{t}(\boldeta))_{i,j}|^{4 q_{1}}\Bigr)=O(1),
\end{equation}
and
%
%
\begin{equation}\label{equB.13}
\max_{1 \leq t \leq n}
\mathrm{E}\Bigl(\sup_{\boldeta\in B_{\delta_{1}}
(\boldeta_{0})}
\|\nabla\varepsilon_{t}(\boldeta)\|^{4 q_{1}}\Bigr)=O(1).
\end{equation}
These equations, however, can be verified based on (\ref{equ3.15}),
(\ref{equ3.28}), (\ref{equB.1}), (\ref{equB.3}), (\ref{equB.4}) and an
argument similar to (\ref{equ3.10})
of~\cite{14}. The details are thus omitted here.
\end{pf*}
\end{appendix}

\section*{Acknowledgments}
We would like to thank the Editor, an Associate Editor and
an anonymous referee for helpful comments and suggestions, which lead
to an improved
version of this paper.


%

%
\printaddresses

\end{document}